\documentclass[11pt,bezier,amstex]{article}  
\usepackage[textsize=small]{todonotes}       
\usepackage{color}              
\usepackage{graphicx}
\usepackage[]{amsmath}
\usepackage{amssymb}
\usepackage{hyperref,bbm}
\definecolor{MyDarkBlue}{rgb}{0,0.08,0.50}
\definecolor{BrickRed}{rgb}{0.65,0.08,0}

\hypersetup{
colorlinks=true,       
    linkcolor=MyDarkBlue,          
    citecolor=BrickRed,        
    filecolor=red,      
    urlcolor=cyan           
}

\newtheorem{Lemma}{Lemma}[section]
\newtheorem{lemma}{Lemma}[section]
\newtheorem{Proposition}[Lemma]{Proposition}

\newtheorem{Theorem}[Lemma]{Theorem}

\newtheorem{Remark}[Lemma]{Remark}

\newtheorem{Construction}[Lemma]{Construction}

\newtheorem{Definition}[Lemma]{Definition}
\newtheorem{Condition}[Lemma]{Condition}

\newcommand{\CA}{\mathcal{A}}

\newcommand{\prob}{\mathbb{P}}
\newcommand{\PP}{\mathcal{P}}

\newcommand{\EE}{\mathcal{E}}
\newcommand{\HH}{\mathcal{H}}
\newcommand{\CC}{\mathcal{C}}
\newcommand{\DD}{\mathcal{D}}
\newcommand{\FF}{\mathcal{F}}
\newcommand{\GG}{\mathcal{G}}
\newcommand{\BB}{\mathcal{B}}
\newcommand{\II}{\mathcal{I}}

\newcommand{\WW}{\mathcal{W}}
\newcommand{\WWstar}{\mathcal{W}^{\star}}
\newcommand{\VV}{\mathcal{V}}
\newcommand{\eps}{\varepsilon}

\newcommand{\var}{{\rm var}\ }

\newcommand{\qed}{\ \ \rule{1ex}{1ex}}

\newcommand{\set}[1]{\left\{#1\right\}}

\newcommand{\Rbold}{{\mathbb{R}}}

\newcommand{\Nbold}{\mathbb{N}}
\newcommand{\R}{{\mathbb R}}
\newcommand{\Zbold}{{\mathbb{Z}}}

\newcommand{\ind}[2]{1_{(e \in \pi(#1,#2))}}

\newcommand{\expec}{\mathbb{E}}

\def\ind{{\rm 1\hspace{-0.90ex}1}}

\newcommand\1{\mathbbm{1}}
\newcommand{\indic}[1]{\1_{\{#1\}}}
\newcommand{\indicwo}[1]{\1_{#1}}
\newcommand{\e}{{\mathrm e}}
\newcommand{\bfd}{\boldsymbol{d}}
\newcommand{\alphan}{\alpha_n}
\newcommand{\gamman}{\gamma_n}
\newcommand{\tn}{\bar{t}_n}
\newcommand{\Kn}{K_n}

\newcommand{\BPbar}{\overline{\BP}}

\newcommand{\eqn}[1]{\begin{equation} #1 \end{equation}}
\newcommand{\eqan}[1]{\begin{align} #1 \end{align}}
\newcommand{\sss}{\scriptscriptstyle}
\newcommand{\Var}{{\rm Var}}
\newcommand{\SWG}{{\sf SWG}}
\newcommand{\BP}{{\sf BP}}
\newcommand{\BPstar}{{\sf BP}^{\star}}
\newcommand{\Alive}{{\sf A}}

\newcommand{\EXP}{{\rm Exp}}

\newcommand{\op}{o_{\sss \prob}}
\newcommand{\Op}{O_{\sss \prob}}
\newcommand {\convp}{\stackrel{\sss {\mathbb P}}{\longrightarrow}}
\newcommand {\convas}{\stackrel{\sss a.s.}{\longrightarrow}}
\newcommand {\convd}{\stackrel{d}{\longrightarrow}}
\newcommand{\proof} {\noindent {\bf Proof}. \hspace{2mm}}

\newcommand {\vep}{\varepsilon}
\newcommand{\nn}{\nonumber}
\newcommand{\MM}{\mathcal{M}}

\numberwithin{equation}{section}
\newcommand{\Ver}{U}
\newcommand{\tildeX}{\widetilde{X}}
\newcommand{\tildeS}{\widetilde{S}}

\newcommand{\GF}{G}
\newcommand{\FR}{F_{\sss R}}
\newcommand{\fR}{f_{\sss R}}
\newcommand{\Hn}{H_n}
\newcommand{\Wn}{L_n}
\newcommand{\tildeHn}{\bar{H}_n}
\newcommand{\tildeWn}{\bar{L}_n}
\newcommand{\Istar}{I^{\star}}
\newcommand{\Fw}{F_{\sss W}}
\newcommand{\Fnw}{F_{n,w}}
\newcommand{\CMnd}{{\rm CM}_n(\bfd)}
\newcommand{\UGnd}{{\rm UG}_n(\bfd)}
\newcommand{\whp}{whp{}}
\renewcommand{\Re}{{\rm Re}}
\newcommand{\Bn}{B_n}
\newcommand{\Cn}{C_n}

\newcommand{\ra}{\rightarrow}
\newcommand{\symdiff}{\triangle}
\newcommand{\mnup}{\overline{m}_n}
\newcommand{\mnlow}{\underline{m}_n}
\newcommand{\YBP}{Y^{\sss{\rm (BP)}}}
\newcommand{\YSWG}{Y^{\sss{\rm (SWG)}}}
\newcommand{\Mis}{{\sf MIS}}
\newcommand{\Tcol}{T^{\sss{\rm (col)}}}
\newcommand{\barTcol}{\bar{T}^{\sss{\rm (col)}}}

\newcommand{\Fstar}{F^{\star}}

\newcommand{\TVD}{d_{\sss \rm TV}}

\newcommand{\ch}[1]{{#1}}
\newcommand{\Scal}{{\mathcal S}}
\newcommand{\Mcal}{{\mathcal M}}
\newcommand{\PiFR}{\Pi^{\sss{\rm (FR)}}}
\newcommand{\PiSR}{\Pi^{\sss{\rm (SR)}}}

\begin{document}
    \author{Shankar Bhamidi
    \thanks{
    Department of Statistics and Operations Research
    The University of North Carolina
    304 Hanes Hall
    Chapel Hill, NC 27510. email: {\tt bhamidi@email.unc.edu}}
    \and
    Remco van der Hofstad
    \thanks{Department of Mathematics and
    Computer Science, Eindhoven University of Technology, P.O.\ Box 513,
    5600 MB Eindhoven, The Netherlands. email: {\tt
    rhofstad@win.tue.nl}}
    \and
    Gerard Hooghiemstra
    \thanks{DIAM, Delft University of Technology, Mekelweg 4, 2628 CD Delft, The
     Netherlands, email: {\tt g.hooghiemstra@ewi.tudelft.nl}}}

\title{Universality for first passage percolation on sparse random graphs}

\maketitle

\begin{abstract}
We consider first passage percolation on sparse random graphs with prescribed
degree distributions and general independent and identically distributed
edge weights assumed to have a density. Assuming that the degree distribution satisfies a uniform $X^2\log{X}$-condition, 
we analyze the asymptotic distribution for the minimal weight
path between \ch{a pair of} typical vertices, as well the number of edges on this path
or hopcount.

The hopcount satisfies a central limit theorem where the \ch{norming} constants are
expressible in terms of the \ch{parameters}  of an associated continuous-time branching
process.
\ch{ Centered
by a multiple of $\log{n}$, where the constant is the inverse of
the Malthusian rate of growth of the associated branching process, the minimal weight
converges in distribution.} The limiting random variable equals
the sum of the logarithms of the martingale
limits of the branching processes that measure the relative growth of
neighborhoods about the two vertices, and a Gumbel random variable,
\ch{and thus shows a remarkably universal behavior.}
The \ch{proofs rely} on a refined coupling between the shortest path problems
on these graphs and continuous-time branching processes, and on a
Poisson point process limit for the potential closing edges of shortest-weight paths between
the source and destination.

The results extend to a host of related random graph models, ranging from
random $r$-regular graphs, inhomogeneous random graphs and
uniform random graphs with a prescribed degree sequence.

\end{abstract}


\noindent
{\bf MSC2000 subject classification:}
60C05, 05C80, 90B15.

\section{Introduction and results}
\label{sec-int}

\subsection{Motivation}
\label{sec-motiv}
First passage percolation is an important topic in modern probability,
due to the inherent applications in a number of fields such as disordered systems in
statistical physics, and since it arises as a building block in the analysis of many more complicated interacting particle systems such as the contact process, other epidemic models and the voter model.

Let us start by describing the basic model. Let $\GG$ be a finite simply
connected graph on $n$ vertices (for example the box $[-N, N]^d$ in the
$\Zbold^d$ lattice, so that $n=(2N+1)^d$). We assign a random edge
\emph{weight or length} independently and identically distributed (i.i.d.)
\ch{to} each of the edges. Due to the random weights, this is an example
of a disordered system entrusted with carrying flow between vertices
in the graph. Fix two vertices in $\GG$. Two functionals of interest
are the minimal weight $\Wn$ of a path between the two vertices and
the number of edges $\Hn$ on the minimal path, often referred
to as the {\it hopcount}. \ch{ We assume that the common distribution of the edge weights is
continuous}, so that the optimal paths are unique and one
can talk about objects such as the number of edges on {\it the} optimal path.

This model has been \ch{intensively} studied, largely in the context of the integer lattice
$[-N,N]^d$ (see e.g.\ \cite{kesten1986aspects, hammersley1966first, 
howard2004models, smythe1978first}). For the power of this model
to analyze more complicated interacting particle systems, see 
\cite{liggett2004interacting} and \cite{durrett1988lecture} and 
the references therein. Due to the interest in
complex networks, such as social networks or the Internet, recently,
this model has attracted attention on general random graph models.
Indeed, stimulated by the availability of an enormous amount of empirical
data on real-world networks, the last decade has witnessed the formulation
and development of many new mathematical graph models for real-world
networks. These models are used to study various dynamics, such as
models of epidemics or random walks to search through the network
(see e.g.\ \cite{albert2002statistical,newman2003j}).

In the modern context, first passage percolation problems take on an
added significance. Many real-world networks (such as the Internet at
the router level or various road and rail networks) are entrusted with
carrying flow between various parts of the network. These networks have
both a graph theoretic structure as well as weights on edges, representing
for example congestion. In the applied setting
understanding properties of both the hopcount and the optimal weight
are crucial, since whilst routing is done via least weight paths, the
actual time delay experienced by users scales like the hopcount
(the number of ``hops'' a message has to perform in getting from
the source to the destination).  \ch{Simulation-based} studies
(see e.g., \cite{braunstein2003optimal}) suggest that random edge
weights have a marked effect on the geometry of the network.
This has been rigorously established in various works
\cite{AmiDraLel11, BhaHofHoo09a, BhaHofHoo09b, BhaHofHoo10},
in the specific situation of \emph{exponential} edge weights.

In this paper, we study the behavior of the hopcount and minimal
weight in the setting of  random graphs with with finite variance degrees and  \emph{general}
continuous edge weights. Since in many applications, the distribution
of edge weights is unknown, the assumption of general weights is 
highly relevant. From a mathematical point of view, working with general
instead of exponential edge weights implies that our exploration process
is \emph{non-Markovian}. This is the first paper that studies first passage percolation
on random graph models in this general setting.
Further, due to the
\ch{choices} of degree distribution, our results immediately carry over to 
various other random graph models, such as rank-1 inhomogeneous random graphs
as introduced in \cite{BolJanRio07}.

{\bf Organization of this section.}
We start by introducing the configuration model
in Section \ref{sec:model-form}, where we also state our main result
Theorem \ref{thm-main-first}.
In Section \ref{sec:math-def}, we discuss a continuous-time branching
process approximation, which is necessary in order to be able to
identify the limiting variables in Theorem \ref{thm-main-first}.
In Section \ref{sec-rel-models},
we extend our results to related random graph models, and in
Section \ref{sec-exam} we study some examples that allow
us to relate our results to results in the literature.
We close with Section \ref{sec-disc}
where we present a discussion of our results and some open problems.

Throughout this paper, we make use of the following standard notation.
We let $\stackrel{a.s.}{\longrightarrow}$ denote convergence almost surely,  $\stackrel{\mbox{L}^1}{\longrightarrow}$ denote
convergence in means, $\convd$ denote convergence in distribution, and 
$\convp$ convergence in probability. For a sequence of random variables
$(X_n)_{n\geq 1}$, we write $X_n=\Op(b_n)$ when $|X_n|/b_n$
is a tight sequence of random variables,
and $X_n=\op(b_n)$ when $|X_n|/b_n\convp0$ as $n\rightarrow\infty$.
To denote that the random variable $D$ has distribution function $F$, we write
$D\sim F$.
For  non-negative functions $n\mapsto f(n)$, $n\mapsto g(n)$
we write $f(n)=O(g(n))$ when $|f(n)|/g(n)$ is uniformly bounded, and
$f(n)=o(g(n))$ when $\lim_{n\rightarrow \infty} f(n)/g(n)=0$.
Furthermore, we write $f(n)=\Theta(g(n))$ if $f(n)=O(g(n))$ and $g(n)=O(f(n))$.
Finally, we write that a sequence of events $(\EE_n)_{n\geq 1}$
occurs \emph{with high probability} (whp) when $\prob(\EE_n)\rightarrow 1$.

\subsection{Configuration model and main result}
\label{sec:model-form}
We work with the configuration model on $n$ vertices $[n]:=\{1,2,\ldots, n\}$.
Let us first describe the random graph model.

\paragraph{Configuration model.}
We are interested in constructing a random graph on $n$ vertices
with prescribed degrees. Given  a
{\it degree sequence}, namely a sequence of $n$ positive integers
$\bfd = (d_1,d_2,\ldots, d_n)$
with $\sum_{i\in [n]} d_i$ assumed to be even, the configuration model (CM) on $n$
vertices with degree sequence $\bfd$ is constructed as follows:

Start with $n$ vertices and $d_i$ half-edges 
adjacent to vertex $i$. The graph is constructed by randomly pairing
each half-edge to some other half-edge to form edges. Let
    \eqn{
    \ell_n = \sum_{i\in [n]} d_i,
    }
denote the total degree. 
Number the half-edges from $1$ to $\ell_n$ in some arbitrary order. 
Then, at each step,
two half-edges which are not already paired are chosen uniformly at random among
all the unpaired or {\it free} half-edges and \ch{have been} paired to form a single edge
in the graph. These half-edges are no longer free and removed from the list of
free half-edges. We continue with this procedure of choosing and pairing two half-edges
until all the half-edges are paired. Observe that the order in which we choose the
half-edges does not matter. Although self-loops may occur,
these become rare as $n\to\infty$ (see e.g.\ \cite{Boll01} or \cite{Jans06b}
for more precise results in this direction). We denote the resulting graph by
$\CMnd$, its vertex set by $[n]$ and its edge set by $\EE_n$.

\paragraph{Regularity of vertex degrees.}
Above, we have described the construction of
the CM when the degree sequence is given. Here,
we shall specify how we construct the actual degree sequence $\bfd$. We start by formulating conditions on $\bfd$.
We denote the degree of a uniformly chosen vertex $V$ in $[n]$ by $D_n=d_{V}$.
The random variable $D_n$ has distribution function $F_n$ given by
    \eqn{
    \label{def-Fn-CM}
    F_n(x)=\frac{1}{n} \sum_{j\in [n]} \indic{d_j\leq x}.
    }
Write $\log(x)_{+}=\log(x)$ for $x\geq 1$ and $\log(x)_{+}=0$ for $x\leq 1$. 
We assume that the vertex degrees satisfy the following
\emph{regularity conditions:}

\begin{Condition}[Regularity conditions for vertex degrees]
\label{cond-degrees-regcond}
~\\
{\bf (a) Weak convergence of vertex \ch{degree}.}\\
There exists a distribution function $F$ on ${\mathbb N}$ such that
    \eqn{
    \label{Dn-weak-conv}
    D_n\convd D,
    }
where $D_n$ and $D$ have distribution functions $F_n$ and $F$, respectively.\\
Equivalently, for any \ch{continuity point $x$ of $F$},
    \eqn{
    \label{conv-Fn-CM}
    \lim_{n\rightarrow \infty} F_n(x)=F(x).
    }
{\bf (b) Convergence of second moment vertex degrees.}
    \eqn{
    \label{conv-sec-mom-Wn}
    \lim_{n\rightarrow \infty}\expec[D_n^2]=\expec[D^2],
    }
where $D_n$ and $D$ have distribution functions $F_n$ and $F$, respectively,
and \ch{we assume that}
	\eqn{
	\label{nu-super-crit}
	\nu=\expec[D(D-1)]/\expec[D]>1.
	}
{\bf (c) Uniform $X^2\log{X}$ condition.} For every $\Kn\rightarrow \infty$,
    \eqn{
    \limsup_{n\rightarrow \infty}\expec[D_n^2\log{(D_n/\Kn)}_+]=0.
    }
\end{Condition}
\ch{The degree of a vertex chosen uniformly at random has distribution $D_n$ as given in Condition \ref{cond-degrees-regcond}(a). By Condition \ref{cond-degrees-regcond}(c), the degree distribution $D_n$ satisfies a uniform  $X^2\log{X}$ condition. A vertex incident to a half-edge that is chosen uniformly at 
random from all half-edges has the same distribution as the random variable
$D^{\star}_n$ given in \eqref{def_Dnstar}, which is the size-biased version of $D_n$. 
The latter random variable satisfies a uniform  $X\log{X}$ condition if and only if 
$D_n$ satisfies a uniform  $X^2\log{X}$ condition. As explained in more detail in 
Section \ref{sec:math-def} below, $D^{\star}_n$ is closely related to a 
\emph{branching-process approximation} of neighborhoods of a uniform vertex, and
thus Condition \ref{cond-degrees-regcond}(c) implies that this branching process
satisfies a uniform $X\log{X}$ condition. By uniform integrability, 
Condition \ref{cond-degrees-regcond}(c) follows from the assumption that $\lim_{n\rightarrow \infty}\expec[D_n^2\log{(D_n)}_+]=\expec[D^2\log{(D)}_+]$.
}

Note that that Condition \ref{cond-degrees-regcond}(a) 
and (c) imply that $\expec[D_n^i]\to \expec[D^i],\, i=1,2$.
When the degrees are \emph{random} themselves, then
we assume that the above convergence conditions hold
\emph{in probability.} Condition \eqref{nu-super-crit}
is equivalent to a giant component existing in $\CMnd$, see
e.g.\ \cite{JanLuc07,MolRee95,MolRee98}.
We often abbreviate $\mu=\expec[D]$.
Let $F$ be a distribution function of a random variable $D$,
 satisfying that $\expec[D^2\log{(D)}_+]<\infty$. 
We give two canonical examples in which
Condition \ref{cond-degrees-regcond} hold. The first is when
there are precisely $n_k=\lceil nF(k)\rceil-\lceil nF(k-1)\rceil$
vertices having degree $k$. The second is when $(d_i)_{i\in [n]}$
is an i.i.d.\ sequence of random variables with distribution
function $F$ (in the case where $\sum_{i\in [n]} d_i$ is odd, we
increase $d_n$ by 1, this \ch{does not} affect the results).

As we \ch{will} describe in more detail in Section \ref{sec-rel-models},
Condition \ref{cond-degrees-regcond} is such that it allows
to extend our results to a range of other random graph models. 

\paragraph{Edge weights and shortest paths.} Once the graph has
been constructed, we attach edge weight $X_e$ to every edge $e$,
where $(X_e)_{e\in \mathcal{E}_n}$ are  i.i.d.~continuous random variables with density
$g\colon [0,\infty)\to [0,\infty)$ and corresponding distribution function $\GF$.
Pick two vertices $\Ver_1$ and $\Ver_2$ at random and let $\Gamma_{12}$ denote the
set of all paths in $\CMnd$ between these two vertices. For any path $\pi \in \Gamma_{12}$,
the weight of the path is defined as
      \begin{equation}
      c_{\tt wk-dis}(\pi)= \sum_{e\in \pi} X_e.
      \label{eqn:wk-dis}
      \end{equation}
Let
    \eqn{
    \label{Wn-def}
    \Wn=\min_{\pi\in \Gamma_{12}}  c_{\tt wk-dis}(\pi),   }
denote the weight
of the optimal (i.e., minimal weight) path between the
two vertices and let $\Hn$ denote number of edges or the
\emph{hopcount} of this path. If the two vertices are in different
components of the graph, then we let $\Wn, \Hn = \infty$.
Now we are ready to state our main result. Due to the complexity of the various constructs (constants and limiting random variables) arising in the theorem, we defer a complete description of these constructs to the next section.

\begin{Theorem}[Joint convergence hopcount and weight]
\label{thm-main-first}
Consider the configuration model $\CMnd$ with degrees $\bfd$ satisfying
Condition \ref{cond-degrees-regcond}, and with i.i.d.\ edge weights
distributed according to the continuous distribution $\GF$. Then, there exist
constants $\alpha, \gamma, \beta\in (0,\infty)$ and $\alphan,\gamman$
with $\alphan\ra \alpha, \gamman\ra \gamma$, such that the hopcount $\Hn$
and weight $\Wn$ of the optimal path between two uniformly selected
vertices conditioned on being connected, satisfy
    \eqn{
    \label{eq-joint-conv}
    \Big(\frac{H_n - \gamman\log{n}}{\sqrt{\beta\log{n}}},~
    \Wn - \frac{1}{\alphan} \log{n}\Big)\convd (Z,Q),
    }
as $n\to\infty$, where $Z$ and $Q$ are independent and $Z$ has a standard normal distribution,
while $Q$ has a continuous distribution.
\end{Theorem}

In Remark \ref{rem-asymp-mean} below, we will state conditions that imply
that we can replace $\alphan$ and $\gamman$ by their limits $\alpha$ and $\gamma$,
respectively.
Theorem \ref{thm-main-first} shows a remarkable kind of universality.
For the configuration model with finite-variance degrees satisfying
Condition \ref{cond-degrees-regcond}, the hopcount always satisfies
a central limit theorem with mean and variance proportional to $\log{n}$.
Also, the weight of the shortest weight path between two
uniformly chosen vertices always is of order $\log{n}$, and the
fluctuations converge in distribution. We will see that even
the limit distribution $Q$ of $\Wn$ has a large degree of
universality. In order to do this, as well as to define the
parameters $\alpha, \alphan,\beta, \gamma, \gamman$, we first need to describe
a continuous-time branching process approximation for
the flow on the configuration model with i.i.d.\ edge weights.

\subsection{Continuous-time branching processes}
\label{sec:math-def}
Before stating our results we recall a standard model of continuous-time
branching process (CTBP), the splitting process or the Bellman-Harris process as well as various associated processes.

Define the size-biased distribution $\Fstar$ of the random variable $D$ with distribution function  $F$ by
    \begin{equation}
    \label{eqn:size-biased}
    \Fstar(x)= \expec[D\indic{D\leq x}]/\expec[D], \quad x\in\R.
    \end{equation}
Now let $(\BPstar(t))_{t\geq 0}$ denote the following CTBP:\\
(a) At time $t=0$, we start with one individual which we shall refer to as the original ancestor or the root of the branching process. Generate $D^{\star}$ having the size-biased distribution $\Fstar$ given in \eqref{eqn:size-biased}.  This individual immediately dies giving rise to $D^{\star}-1$ children.
\\(b) Each new individual $v$ in the branching process lives for a random amount of time which has distribution $\GF$, i.e., the edge weight distribution, and then dies. At the time of death again the individual  gives birth to $D_v^{\star}-1$ children, where $D_v^{\star}\sim
\Fstar$. Lifetimes and number of offspring across individuals are independent.

Note that in the above construction, by Condition \ref{cond-degrees-regcond}(b), if we let $X_v= D_v^{\star}-1$ be the number of children of an individual then the expected number of children satisfies
    \begin{equation}\label{eqn:exp-child}
    \expec[X_v] = \expec[D^{\star}_v-1] = \nu > 1,
    \end{equation}
Further, by
Condition \ref{cond-degrees-regcond}(c), for $D^{\star}\sim \Fstar$,
    \begin{equation}\label{eqn:xlogx}
    \expec[D^{\star}\log(D^{\star})_+] < \infty.
    \end{equation}

The CTBP defined above is a {\it splitting} process, with lifetime distribution $\GF$ and offspring
distribution $D^{\star}-1$. Denote by $N(t)$ the number of offspring of an individual at time $t$. Then
$\mu(t)=\expec[N(t)]=\nu\cdot \GF(t)$, where $\nu=\expec[D^{\star}-1]$. Furthermore, the Laplace-Stieltjes transform ${\hat \mu}$
is defined by
    $$
    {\hat \mu}(s)=\nu \int_0^\infty \e^{-st}\,d\GF(t).
     $$
The Malthusian parameter $\alpha$ of the branching process $\BPstar(\cdot)$ is the unique solution of the equation
    \begin{equation}
    \label{eqn:malthusian-param}
    {\hat \mu}(\alpha)=\nu \int_0^\infty  \e^{-\alpha t} d\GF(t)=1.
    \end{equation}
Since $\nu>1$, we obtain
that $\alpha\in (0,\infty)$. We also let $\alphan$ be the solution to
\eqref{eqn:malthusian-param} with $\nu$ replaced with $\nu_n=\expec[D_n(D_n-1)]/\expec[D_n].$
Clearly, $\alphan\ra \alpha$, when Condition \ref{cond-degrees-regcond} holds, 
and $|\alphan-\alpha|=O(|\nu_n-\nu|)$.

Standard theory (see e.g., \cite{athreya-ney-book, jagers1975branching, jagers1984growth}) implies that under our assumptions of the model, namely \eqref{eqn:exp-child} and \eqref{eqn:xlogx}, there exists a random variable $\WWstar$ such that
    \begin{equation}
    \e^{-\alpha t} |\BPstar(t)|\stackrel{a.s., \mbox{L}^1}{\longrightarrow} \WWstar.
    \label{eqn:bp-asymp}
    \end{equation}
Here the limiting random variable $\WWstar$
satisfies $\WWstar>0$ a.s.\ on the event of non-extinction of the branching
process and is zero otherwise. Thus $\alpha$ measures the true rate of 
exponential growth of the branching process.

Define the distribution function ${\bar \GF}$, which is often referred to as the
\emph{stable-age distribution}, by
    \begin{equation}
    {\bar \GF}(x) = \nu \int_0^x  \e^{-\alpha y} d\GF(y),
    \label{eqn:stable-age-dist}
    \end{equation}
where we recall that $\alpha$ is the Malthusian rate of growth parameter.
Let $\bar{\nu}$ be the mean and $\bar{\sigma}^2$ the variance of
${\bar \GF}$. Then  $\bar{\nu},\bar{\sigma}^2\in (0,\infty)$, since $\alpha>0$.
We also define ${\bar \GF}_n$ to be the distribution function ${\bar \GF}_n$
in \eqref{eqn:stable-age-dist} with $\nu$ and
$\alpha$ replaced with $\nu_n$ and $\alphan$, and we let $\bar{\nu}_n$ and
$\bar{\sigma}^2_n$ be its mean and variance.

We need a small variation of the above branching process where the root of
the branching process dies immediately giving birth to a $D$ number of
children where $D$ has distribution $F$. The details for every other
individual in this branching process remain unchanged from the original
description, namely each individual survives for a random amount of time
with distribution $\GF$ giving rise to a $D^{\star}-1$ number of children
where $D^{\star}\sim \Fstar$, the size-biased distribution function $\Fstar$.
Writing $\ch{|\BP(t)|}$ for the number of alive individuals at time $t$, it
is easy to see here as well that
    \eqn{
	\label{conv-as-BP}
    \e^{-\alpha t}|\BP(t)| \stackrel{a.s., \mbox{L}^1}{\longrightarrow} \tilde{\WW}.
    }
Here, \ch{$\tilde{\WW}$} satisfies the stochastic equation,
	  \[\tilde{\WW}= \sum_{i=1}^D \WW^{\star,\sss(i)} \e^{-\alpha \xi_i},
    \]
where $D\sim F$, and $\tilde{\WW}^{\star,\sss(i)}$ are i.i.d.\ with the 
distribution of the limiting random variable in \eqref{eqn:bp-asymp}, and \ch{$\xi_i$} 
are i.i.d.~with distribution $\GF$.
Let
    \begin{equation}
    \label{eqn:tilde-w}
    \WW =\tilde{\WW}|\tilde{\WW}>0.
    \end{equation}
To simplify notation in the sequel, we will use $(\BP(t))_{t\geq 0}$ to denote a 
CTBP with the root having offspring either $D$ or $D^{\star}-1$, which will
be clear from the context.

We are now in a position to identify the limiting random variable $Q$
as well as the parameters $\alpha, \beta, \alphan, \gamma_n, \gamma$:

\begin{Theorem}[Identification of the limiting variables]
\label{thm-main}
The parameters $\alpha, \alphan,\beta,\gamman$ in Theorem \ref{thm-main-first}
satisfy that $\alpha$ is the Malthusian rate of growth defined in
\eqref{eqn:malthusian-param} and $\alphan$ is the solution 
to \eqref{eqn:malthusian-param} with $\nu_n$ replacing $\nu$, while
    \eqn{
    \gamman=\frac{1}{\alphan \bar{\nu}_n},
    \qquad
    \beta=\frac{\bar{\sigma}^2}{\bar{\nu}^3\alpha}.
    }
Further, $Q$ can be identified as:
    \eqn{
    \label{lim-var}
    Q=\frac{1}{\alpha}\left(-\log{\WW^{\sss(1)}}-\log{\WW^{\sss(2)}}-\Lambda+c\right),
    }
where $\prob (\Lambda\le x)=\e^{-\e^{-x}}$, so that $\Lambda$ is a standard Gumbel random variable,
$\WW^{\sss(1)}, \WW^{\sss(2)}$ are two independent copies of the variable $\WW$
in \eqref{eqn:tilde-w}, also independent from $\Lambda$, and $c$ is the constant
    \eqn{
    c=\log(\mu(\nu-1)^2/(\nu \alpha \bar{\nu})).
    }
\end{Theorem}

\begin{Remark}[Asymptotic mean]
\label{rem-asymp-mean}
We can replace $\alphan$ and $\gamman$ by their limits $\alpha$ and \ch{$\gamma=1/(\alpha \bar{\nu})$} in 
\eqref{eq-joint-conv} precisely when $\gamman=\gamma+o(1/\sqrt{\log{n}})$ and
$\alphan=\alpha+o(1/\log{n})$. Since $|\alphan-\alpha|=O(|\nu_n-\nu|)$,
$|\bar{\nu}_n-\bar{\nu}|=O(|\nu_n-\nu|)$,
these conditions are equivalent to $\nu_n=\nu+o(1/\sqrt{\log{n}})$ and 
$\nu_n=\nu+o(1/\log{n})$, respectively.
\end{Remark}

Theorem \ref{thm-main} implies that also the random variable $Q$
is remarkably universal, in the sense that it always
\ch{involves} two martingale limit variables
corresponding to the flow problem, \ch{and} a Gumbel
distribution. While such results were known for the exponential
distribution (see e.g., \cite{BhaHofHoo09a}), this is the first
time that FPP on random graphs with general edge weights is studied.

Let $\Wn(i)$ denote the weight of the $i^{\rm th}$ shortest path,
so that $\Wn=\Wn(1)$, and let $\Hn(i)$ denote its length.
Further let $\tildeHn(i)$ and $\tildeWn(i)$ denote the re-centered 
and normalized quantities as in Theorem \ref{thm-main-first}. \ch{The same proof for the optimal path easily extends to prove asymptotic results for the joint distribution of the weights and hopcount of these ranked paths. To keep the study to a manageable length, we shall skip a proof of this easy extension.  }

\begin{Theorem}[Multiple paths]
\label{thm-main-paths}
Under the conditions of Theorem \ref{thm-main-first},
for every $m\geq 1$,
    \eqn{
    \label{eq-joint-conv2}
    ((\tildeHn(i),\tildeWn(i)))_{i\in[m]}\convd ((Z_i,Q_i))_{i\in [m]},
    }
as $n\to\infty$, where for $i\ge 1$, $Z_i$ and $Q_i$ are independent and $Z_i$ has a standard normal distribution,
while 
	\eqn{
	Q_i=\frac{1}{\alpha}\left(-\log{\WW^{\sss(1)}}-\log{\WW^{\sss(2)}}-\Lambda_i+c\right),
	}
\ch{where $(\Lambda_i)_{i\geq 1}$ are the ordered points of an inhomogeneous Poisson point process with intensity $\lambda(t)=\e^{t}$}. 
\end{Theorem}

\subsection{Related random graph models}
\label{sec-rel-models}
\paragraph{Uniform random graphs with a prescribed degree sequence.}
We call a graph \emph{simple} when it contains no self-loops nor multiple
edges. It is well known that the CM conditioned on
being simple is a \emph{uniform} random graph with the same degrees.
As a result, our \ch{theorems} extend to this setting:

\begin{Theorem}[Extension to uniform random graphs with prescribed degrees]
\label{thm-unif-RGs}
Under\\
 the conditions of Theorem \ref{thm-main-first},
the results in Theorems \ref{thm-main-first} and the identification
of the limiting variables in Theorem \ref{thm-main} apply to
uniform random graphs with prescribed degree sequence $\bfd$ satisfying
Condition \ref{cond-degrees-regcond}.
\end{Theorem}
The proof of Theorem \ref{thm-unif-RGs} follows rather directly from
that of Theorems \ref{thm-main-first}-\ref{thm-main}, by conditioning
on simplicity. By \cite{Boll01} or \cite{Jans06b}, under Condition \ref{cond-degrees-regcond},
    \eqn{
    \label{prob-simple}
    \lim_{n\rightarrow \infty} \prob(\CMnd \text{ simple})=\e^{-\nu/2-\nu^2/4}.
    }
This proof of \eqref{prob-simple} follows by a Poisson approximation on the number of self-loops and
the number of multiple edges, which are proved to converge to two
independent Poisson random variables with means $\nu/2$ and $\nu^2/4$ respectively.
We can interpret the probability in \eqref{prob-simple} as the probability
that both these Poisson variables are equal to zero, which is
equivalent to the graph being simple. Now the proof of the main theorem reveals that in order to find the minimal
weight path between vertices $\Ver_1,\Ver_2$, we only need to investigate
of order $\sqrt{n}$ edges. Therefore, the event of simplicity
of the configuration model will be mainly determined by
the \emph{uninspected} edges, and is therefore asymptotically
independent of $(\Hn,\Wn)$. This explains Theorem \ref{thm-unif-RGs}.
We give a full proof of Theorem \ref{thm-unif-RGs} in Section \ref{sec-pfs-relat-RGs}.

\paragraph{Rank-1 inhomogeneous random graphs.}
Fix a sequence of positive weights $(w_i)_{i\in [n]}$.
We shall assume that there exists a distribution function $\Fw$ on $\Rbold^+$ such that
    \eqn{
    \Fnw(x)= \frac{1}{n}\sum_{i\in[n]} \ind_{\{w_i\leq x \}}\to \Fw(x),
    \label{eqn:cdf-wconv}
    }
for each point of continuity of $\Fw$. Here, $\ind_A$ denotes the indicator of the
set $A$. Let $W_n$ denote the weight of a uniformly chosen vertex in $[n]$, i.e.,
$W_n=w_{V}$, where $V\in [n]$ is chosen uniformly. Let $W_n^{\star}$ denote the
size-biased version of $W_n$, i.e.,
    \eqn{
    \prob(W_n^{\star}\leq x)=\expec[W_n \indic{W_n\leq x}]/\expec[W_n].
    }


Now given these weights, we construct
a random graph by attaching an edge between vertex $i$ and $j$
with probability
    \begin{equation}
    \label{pij-def}
    p_{ij} = 1-\e^{-w_i w_j/\ell_n},
    \end{equation}
where, \ch{with some abuse of notation},
    \eqn{
    \ell_n=\sum_{i\in [n]} w_i,
    }
is the sum of the vertex weights, and the status of different
edges are independent. Let $\nu=\expec[W^2]/\expec[W]$. We always assume $\nu>1$ as this
is necessary and sufficient for the existence of a giant component
(see, e.g., \cite{BolJanRio07}). Note that letting $w_i=\lambda$,
we immediately get the Erd\H{o}s-R\'enyi random graph with edge
connection probability $1-\e^{-\lambda/n}$. Thus, this model is a
natural generalization of the classical random graph model.
Related models are the \emph{generalized random graph}
introduced by Britton, Deijfen and Martin-L\"of in
\cite{BriDeiMar-Lof05}, for which
    \eqn{
    \label{pij-GRG}
    p_{ij} = \frac{w_i w_j}{\ell_n+w_iw_j},
    }
and the \emph{random graph with given prescribed degrees}
or Chung-Lu model, where instead
    \eqn{
    \label{pij-CL}
    p_{ij}=\max(w_iw_j/\ell_n, 1),
    }
and which has been studied intensively by Chung and Lu (see
\cite{ChuLu02a,ChuLu02b,ChuLu03, ChuLu06c, ChuLu06}).
\ch{Let $W$ denote a random variable with distribution $\Fw$.}
By Janson \cite{Jans08a}, when $W_n\convd W$ and
$\expec[W_n^2]\rightarrow \expec[W^2]$, the three random graph models defined above
are \emph{asymptotically equivalent}, meaning that all events
have asymptotically equal probabilities.
By \cite{BolJanRio07}, with $N_k(n)$ denoting
the number of vertices with degree $k$,
    \eqn{
    \label{degree-conv-IRG}
    N_k(n)/n \convp \expec\left[\e^{-W} \frac{W^k}{k!}\right].
    }
This proves that $D_n\convd D$,
where $D$ has the mixed Poisson distribution given in \eqref{degree-conv-IRG}, i.e.,
for this model
	\eqn{
	\label{F-IRG}
	F(x)=\sum_{k\leq x} \expec\left[\e^{-W} \frac{W^k}{k!}\right].
	}
This is formulated in the following theorem:

\begin{Theorem}[Extension to rank-1 inhomogeneous random graphs]
\label{thm-IRGs}
For rank-1 inhomogeneous random graphs with edge probabilities in \eqref{pij-def},
\eqref{pij-GRG} or \eqref{pij-CL}, where the weight of a uniform vertex
$W_n$ satisfies that $W_n\convd W$, $\expec[W_n]\rightarrow \expec[W],
\expec[W_n^2]\rightarrow \expec[W^2]$ and, for every $\Kn\ra \infty$,
    \[
    \ch{\expec[W^2_n\log{(W_n/\Kn)}_+]=o(1),}
    \]
the results in Theorems \ref{thm-main-first} and \ref{thm-main} hold
with limiting degree distribution $F$ in \eqref{F-IRG}.
\end{Theorem}

Theorem \ref{thm-IRGs} can be understood as follows. By \cite{BriDeiMar-Lof05},
the generalized random graph \emph{conditioned on its degree sequence}
is a uniform random graph with the same degree sequence. Therefore,
Theorem \ref{thm-IRGs} follows from Theorem \ref{thm-unif-RGs} when
the conditions on the degrees in Condition \ref{cond-degrees-regcond}
hold in probability. By \eqref{degree-conv-IRG}, $D_n\convd D$,
where $D$ has the mixed Poisson distribution given in \eqref{degree-conv-IRG}.
Therefore, it suffices to show that $\expec[W_n]\rightarrow \expec[W],
\expec[W_n^2]\rightarrow \expec[W^2]$ and $\lim_{n\ra \infty} \expec[W_n^2\log{(W_n/\Kn)}_+]=0$
imply that the same convergence holds for the degree of a random vertex.
This is proved in Section \ref{sec-pfs-relat-RGs}.

\subsection{Examples}
\label{sec-exam}
In this section, we discuss a few special examples of the edge 
weights that have arisen in a number of different contexts in the 
literature and have been treated via distribution specific techniques.

\paragraph{Exponential weights.}
FPP on random graphs with exponential edge weights have received substantial
attention in the literature (see e.g., \cite{AmiDraLel11, BhaHofHoo09a, BhaHofHoo09b, BhaHofHoo10}).
Let $\GF(x)=1-\e^{-x},$ for $x\geq 0$, denote the distribution function
of an exponential random variable with mean 1. This was one of the first
models to be formulated and analyzed in the context of the integer lattice,
see \cite{morgan1965two} and the complete analysis in \cite{hammersley1966first}.
For exponential weights, the Malthusian rate of growth parameter $\alpha$ satisfies
 \begin{equation}
    \nu \int_0^\infty  \e^{-\alpha x} d\GF(x)=\nu \int_0^\infty  \e^{-\alpha x} \e^{-x}dx=1,
    \end{equation}
so that $\alpha=\nu-1,\alphan=\nu_n-1$. Similarly, one can compute that 
${\bar \GF}_n(x)=1-\e^{-\nu_n x},\, x\geq 0$, so that
    \eqn{
    \bar{\nu}_n= 1/\nu_n, \qquad \bar{\nu}= 1/\nu,\qquad \bar{\sigma}^2= 1/\nu^2.
    }
Using these values in Theorem \ref{thm-main} shows that $\Hn$ converges to a normal distribution, with
asymptotic mean and asymptotic variance both equal to $\frac{\nu}{\nu-1}\log n$.
Finally $c=\log{\Big(\mu(\nu-1)^2/(\alpha \nu \bar{\nu})\Big)}=\log(\mu(\nu-1))$,
which is equal to the constant in \cite[(3.7)]{BhaHofHoo09a}.
\footnote{In \cite[(3.7)]{BhaHofHoo09a}, the Gumbel variable $\Lambda$,
which appears as $\log{M}$ in \cite[(C.19)]{BhaHofHoo09a}, 
where $M$ is an exponential random variable,
should be replaced with $-\Lambda$.}
\footnote{Also in \cite[(5.4)]{BhaHofHoo09a} there is an error in the
precise limiting random variable for FPP on the Erd\H{o}s-R\'enyi 
random graph due to the fact that \cite[(4.16)]{BhaHofHoo09a}  is not correct.}
This result thus generalizes the hopcount result in \cite{BhaHofHoo09a},
where FPP with exponential weights is considered on the configuration model
with i.i.d.\ degrees with $D\geq 2$ a.s. In \cite{BhaHofHoo09a}, also
\emph{infinite variance} degrees are studied, a case that we do not
investigate here. In fact, some of the results we prove here
do \emph{not} extend to this setting, see Section \ref{sec-disc}
for more details. 


\paragraph{Exponential weights plus a large constant.}
We next study what happens when $X_e=1+E_e/k$, where $(E_e)$ are i.i.d.\ exponentials with mean 1, and
$k$ is a large constant. This setting is, apart from a trivial
time-rescaling, identical to the setting where $X_e=k+E_e$.
In this case, one would expect that for large $k$, $\Hn$ is close
to the graph distance between a pair of uniformly chosen vertices in $[n]$,
conditioned to be connected. This problem has
attracted considerable attention, see, in particular, \cite{EskHofHoo06} for the
Norros-Reittu model and \cite{HofHooVan05a} for the CM. In these works, it has been shown that \ch{$(\Hn-\log_\nu{(n)})_{n\ge 1}$ is a tight sequence of random variables}. This suggests (compare with Theorems
\ref{thm-main-first}-\ref{thm-main}) that,
as $k\rightarrow \infty$,
    \eqn{
    \alpha\to \log{\nu},
    \qquad \bar{\nu}\to 1,
    \qquad
    \frac{\bar{\sigma}^2}{\bar{\nu}^3\alpha}\to 0.
    }
We now check this intuitive argument. Indeed,
    \eqn{
    \nu \int_0^\infty  \e^{-\alpha x} d\GF(x)
    =\nu k\int_1^\infty  \e^{-\alpha x} \e^{-k(x-1)}dx
    =\frac{\nu k}{\alpha+k}\e^{-\alpha}=1.
    }
While solving this equation \emph{explicitly} is hard,
it is not too difficult to see that $k\to \infty$ implies that
$\alpha\to \log{\nu}$.

We can compute
the stable-age distribution as
$1+\EXP(k+\alpha)$,
so that $\bar{\nu}=1+\frac{1}{k+\alpha}$, while
$\bar{\sigma}^2=1/(k+\alpha)^2\rightarrow 0$.
Therefore, $\bar{\nu}\sim 1$, which in turn also implies that
$\alpha \bar{\nu}\to \log{\nu}$.
Also,
    \eqn{
    \frac{\bar{\sigma}^2}{\bar{\nu}^3\alpha}\sim k^{-2}(\log{\nu})^{-1}\rightarrow 0.
    }
This shows that the two settings of graph distances and FPP with weights $1+\EXP(1)/k$
match up nicely when $k\rightarrow \infty$.


\paragraph{Weak disorder on random regular graph with large degree.}
\ch{As a third example we} consider the configuration model with fixed degrees $r$, and where each edge is
given an edge weight $E^s,\, s>0$ where $E\sim \EXP(1)$. The parameter $s$ plays the role of inverse temperature in statistical physics with $s\to \infty$ corresponding to the minimal spanning tree with exponential edge weights. This setting has been studied on the complete graph in
\cite{BhaHof10}, and here we make the connection to the results proved there. 

In this case, $\nu=r-1$.
The Malthusian parameter $\alpha$ satisfies (compare \eqref{eqn:malthusian-param}),
    \eqn{
    \label{malpar}
    (r-1)p\int_0^\infty \e^{-\alpha x-x^p}x^{p-1}\,dx=1,
    }
where $p=1/s$. We can not solve \eqref{malpar} explicitly, but when
$r\rightarrow \infty$, we \ch{conclude that} $\alpha\rightarrow \infty$, so that the above equation
is close to
\ch{
    \eqn{
    \label{malpar-large-d}
    r p\int_0^\infty \e^{-\alpha x}x^{p-1}\,dx=r p \alpha^{-p}\Gamma(p)=r\alpha^{-p}\Gamma(p+1)=1,
    }
    }
so that
    \eqn{
    \label{alpha-large-d}
    \alpha=(r\Gamma(1+1/s))^{s}.
    }
The moments ${\bar \nu}$ and ${\bar \sigma^2}$ are then approximately given by
 \ch{   \eqn{
    {\bar \nu}\approx r p\int_0^\infty \e^{-\alpha x}x^{p}\,dx
    =r p \alpha^{-(p+1)}\Gamma(p+1)=p/\alpha,
    }
    }
and\ch{
    \eqn{
    {\bar \sigma^2}\approx r p\int_0^\infty \e^{-\alpha x}x^{p+1}\,dx-{\bar \nu}^2
    =rp \alpha^{-(p+2)} \Gamma(p+2)-(p/\alpha )^2=p/\alpha^2,
    }
    }
where we repeatedly use \eqref{alpha-large-d}.
As a result, we obtain
    \eqn{
    \gamma=\frac{1}{\alpha \bar{\nu}}\approx 1/p=s,
    \qquad
    \beta=\frac{\bar{\sigma}^2}{\bar{\nu}^3\alpha}\approx 1/p^2=s^2.
    }
These results match up nicely with the result on the complete graph,
obtained when \ch{$r=n-1$}, for which \cite{BhaHof10} show that
a central limit theorem holds for $\Hn$ with asymptotic mean $s\log{n}$ and
asymptotic variance $s^2 \log{n}$, while $n^s[\Wn-\frac{1}{\lambda}\log{n}]$
converges in distribution, where $\lambda=\Gamma(1+1/s)^{s}$.


\subsection{Discussion}
\label{sec-disc}
In this section, we discuss our results, possible extensions and open problems.\\
(a) {\bf Universality.} As Theorems \ref{thm-main-first}-\ref{thm-main} show, even second
order asymptotics for the hopcount in the presence of disorder in the
network depends only on the first two moments of the size-biased
offspring distribution \ch{and on the edge-weight distribution, but not on any other property of the network model}. Further, the limit distribution of the minimal weight
between two uniform vertices conditioned on being connected
has a universal shape, even though the martingale
limit of the flow naturally strongly depends both on the graph
topology as well as on the edge weight distribution.\\
(b) {\bf Divergence from the mean-field setting.} One famous model
that has witnessed an enormous amount of interest in probabilistic
combinatorial optimization is the mean-field model, where one starts
with the complete graph on $n$ vertices and then attaches random edge weights
and analyzes optimal path structure. See \cite{janson1999one}
for a study on the effect of exponential edge weights on the
geometry of the graph, \cite{frieze1985value} for study of the
minimal spanning tree and \cite{aldous-steele-obj} for a discussion
of a number of other problems. Branching process methods have been
used to good effect in this setting as well to analyze the effect
of random disorder on the geometry of the graph, see \cite{BhaHof10}
for example and often give good heuristics for what one would expect 
in the sparse random graph setting. However, what the main theorems
imply is that in a number of situations, the mean-field setting
diverges markedly from the random graph setting. For example, when
each edge has \ch{$E^{-s}$} weight where $E$ has an exponential
distribution and $s> 0$,
one can show that in fact the hopcount between typical vertices
converges to a constant \cite{BhaHofHoo09a}, while Theorem
\ref{thm-main} implies that even in this case, for the CM, the hopcount scales as $\log{n}$ and satisfies a CLT.\\
(c) {\bf Infinite variance configuration model.} In \cite{BhaHofHoo09b},
we have also investigated the setting where the degrees are i.i.d.\ with
$\expec[D^2]=\infty$ and with exponential edge weights.
In this case, the result for $\Wn$ is markedly
different, in the sense that \ch{ $\Wn$ converges in distribution without re-centering}.
Further, $\Hn$ satisfies a central limit with asymptotic mean
and  variance equal to \ch{ a multiple of $\log{n}$}. Now, when taking
$X_e\sim 1+\EXP(1)$, by \cite{HofHooZna04a}, there are paths of length
$\log\log{n}$ connecting vertices $\Ver_1$ and $\Ver_2$, conditionally on
$\Ver_1$ and $\Ver_2$ being connected. Since the weight of such a path
is of the same order, we conclude that $\Hn, \Wn=\Theta_{\sss \prob}(\log\log{n})$.
Thus, in such \ch{cases} it is possible that $\Hn$ acts on a different scale,
even though $\Hn\convp \infty$. It would be of interest to investigate
whether $\Hn$ always satisfies a central limit theorem,
and, if so, whether the order of magnitude of its variance is always
equal to that of its mean.\\
(d) {\bf The $X\log{X}$-condition.} By Condition \ref{cond-degrees-regcond}(c),
we assume that the degrees satisfy a second moment condition with an \ch{additional} logarithmic
factor. This is equivalent to the CTBP
satisfying an $X\log{X}$ condition (uniformly in the size $n$ of the graph).
It would be of interest to investigate what happens when this condition fails.
It is well known that the branching process martingale limit is identically equal
to 0 when $\expec[X]<\infty$, but $\expec[X\log{(X)}_+]=\infty$
(see e.g., \cite{athreya-ney-book} or \cite{jagers1975branching},
or \cite{LyoPemPer95}). Therefore, the limit in \eqref{lim-var} does not exist.
This suggests that \eqref{eq-joint-conv} should be replaced with
$\Wn - \frac{2}{\alphan}\log{(t_n)}$, where $t_n$ is such that
$|\BPstar(t_n)|n^{-1/2}$ has a non-trivial limit.\\
(e) {\bf Flooding and diameter.} In \cite{AmiDraLel11}, the \emph{flooding time} and \emph{diameter},
i.e., $\max_{j\in [n]\colon L_{\Ver_1,j}<\infty}L_{\Ver_1,j}$, respecively $\max_{i,j\in[n]\colon L_{i,j}<\infty} L_{i,j}$, 
where $L_{i,j}$ is the minimal
weight between the vertices $i$ and $j$ and $\Ver_1$ is, as before, a randomly selected vertex,
is investigated in the context of the CM with exponential edge weights. It would be 
of interest to investigate the flooding time
for general edge weights. We expect that the exponential distribution is special,
since there the typical weight has the same order of magnitude as the maximum
over the vertices of the minimal edge weight from that vertex. This fact is only true
when the weight distribution has an exponential tail. For example, taking $X_e=E_e^{s}$
for $s>1$, the maximal minimal weight from a vertex is of order $(\log{n})^s$, which
is much larger than the typical distance, which is $\Theta_{\sss \prob}(\log{n})$ due to our main results.
It would be interesting to investigate what the limit of the weight diameter is in 
this simple example.
\\
(f) {\bf Superconcentration and chaos.} In this study we have looked at some global functionals of the optimal path between randomly selected vertices and in particular have shown that the weight of the optimal path satisfies \ch{$L_n / \log{n} = \Op(1)$}. Analogous to various related problems in statistical physics such as random polymers, or FPP on the lattice, this suggests that the optimal path problem satisfies superconcentration. \ch{In particular, it suggests that}  this random combinatorial optimization problem is \emph{chaotic} in the sense that there exists $\eps_n\to 0$ such that refreshing a fraction $\eps_n$ of the edge weights with new random variables with the same distribution would entirely change the actual optimal path, \ch{in the sense that} the new optimal path would be ``almost'' disjoint of the original optimal path, see e.g.\ \cite{chatterjee2008chaos}.  Such questions have also arisen in computer science wherein one is interested in judging the ``importance'' and fair price of various edges in the optimal path; if an edge being deleted causes a large change in the cost of the new optimal path, then that edge is deemed very valuable. These form the basis of various ``truth and auction mechanisms'' in computer science (see e.g.\ \cite{mihail2006certain}, \cite{flaxman2006first}, \cite{bhamidi2008first}). It would be interesting to derive rigorous results in our present context. 
\\
(g) {\bf Pandemics, gossip and other models of diffusion:} First passage percolation models as well as models using FPP as a building \ch{block} have started to play an increasingly central role in the applied probability community in describing the flow of materials, ranging from viral epidemics (\cite{epidemics}), gossip algorithms (\cite{aldous2010knowing}) and more general finite Markov interchange processes (\cite{aldous2012lecture}). Models with more general edge distributions have also arisen in understanding the flow of information and reconstruction of such information networks in sociology and computer science, see e.g.\ \cite{leskovec2009meme},\cite{leskoveccascading} for just some examples in this vast field.

\section{Proof: construction of the flow clusters}
\label{sec-pfs}
We start with some central constructions that lay the ground work for the proofs of the main results. 
We denote by $\Ver_1$ and $\Ver_2$
two randomly selected vertices, conditioned on being connected.
We think of the weights as edge lengths so that they induce a random metric on the graph
$\CMnd$.  For a half-edge $y$, we let $P_y$ denote
the half-edge to which it is paired, i.e., $(y,P_y)$ forms an edge. Further, \ch{ we let $V_y$ be the vertex to which the half-edge $y$  is incident.}

\subsection{Flow clusters from $\Ver_1$ and $\Ver_2$} 
\label{sec-flow-clusters}
To understand the shortest path between these vertices, think of water percolating through the network at rate one, started simultaneously from the two vertices. For any $t\geq 0$, the set of vertices first seen by the flow from \ch{$\Ver_i$} will often referred to the flow cluster or the shortest weight graph of vertex $\Ver_i$. When the two flows collide or create prospective collision edges, then these generate prospective shortest \ch{paths}.

Let us now give a precise mathematical formulation to the above description.  
We grow two flow clusters (i.e.\ two stochastic processes in
continuous time) from $\Ver_1$ and $\Ver_2$, simultaneously. The main ingredients of the two flow clusters,
namely the alive set $\Alive(t)$ will only change at random times $T_0=0<T_1<T_2<\ldots$ and therefore
the definition can be given recursively. At time $t=T_0=0$, the vertices $\Ver_1$ and $\Ver_2$ die instantaneously,
and give rise to $d_{\Ver_1}$ and $d_{\Ver_2}$ children. These children correspond to half-edges incident to
$\Ver_1$ and $\Ver_2$. We start by testing whether any of the half-edges incident to $\Ver_1$ are paired
to one another. If so, then we remove both half-edges from the total set of
$d_{\Ver_1}$ half-edges. We then define $X_0^{\sss(1)}$ the number of unpaired half-edges
after the self-loops incident to $\Ver_1$ are removed. We next continue with the $d_{\Ver_2}$
half-edges incident to $\Ver_2$, and check whether they are paired to one of the $X_0^{\sss(1)}$ 
remaining half-edges incident to $\Ver_1$ or any of the $d_{\Ver_2}$ half-edges incident to $\Ver_2$. 
When such a half-edge is paired to one of the $d_{\Ver_2}$ sibling half-edges, 
a self-loop is formed. When such a half-edge is paired to one of the $X_0^{\sss(1)}$ 
remaining half-edges incident to vertex $\Ver_1$, a so-called 
{\it collision edge} is formed. A collision
possibly yields the path with minimal weight between $\Ver_1$ and
$\Ver_2$. We let $X_0^{\sss(2)}$ denote the number of unpaired half-edges after the
tests for collision edges and cycles have been performed. 
Note that, by construction, each of the $X_0^{\sss(i)}$ half-edges 
incident to the vertices $\Ver_i$, where $i\in \{1,2\}$, are paired to 
\emph{new} vertices, i.e., vertices distinct from $\Ver_1$ and $\Ver_2$.

For the moment we collect the collision edges at time $T_0$, together with the
weights of the connecting edge between $\Ver_1$ and $\Ver_2$, and continue with the description
of the flow clusters. All edges that are not paired to one of the other $d_{\Ver_1}+d_{\Ver_2}-1$ half-edges
incident to either $\Ver_1$ or $\Ver_2$ together form the set $\Alive(0)$, the set of active half-edges at time $0$.
For $y\in \Alive(0)$, we define $I(y)=i$ if the half-edge $y$ is connected to $\Ver_i$, $i=1,2$,
and we define $(R_0(y))_{y\in \Alive(0)}$ as an i.i.d.~sequence of life times having
distribution function $G$.

We denote the set of half-edges at time $t$ by $\Alive(t)$. For $y\in \Alive(t)$, we record
$I(y)$, which is the index $i\in \{1,2\}$ to which \ch{$\Ver_i$} the half-edge is connected, and
we let $H(y)$ denote the height of $y$ to \ch{$\Ver_{I(y)}$}. This height
equals $0$ for $y\in \Alive(0)$. When we introduce new half-edges at $\Alive(t)$ at later times
we will specify the height of these half-edges.
Now define $T_1=\min_{y\in \Alive(0)} R_0(y)$ and denote by $y_0^{\star}$ the half-edge
equal to the argument of this minimum, hence $R_0(y^{\star}_0)=\min_{y\in \Alive(0)} R_0(y)$.
Since life-times have a continuous distribution, $y^{\star}_0$ is a.s.~{\it unique}.
Now set $\Alive(t)=\Alive(0),\,0\le t<T_1 $, i.e., the active set remains
the same during the interval $[0,T_1)$, and define the flow cluster $\SWG(t)$, for $0\leq t <T_1 $, by
	\eqn{
	\label{def-flow}
	\SWG(t)=\{y,I(y),H(y),R_t(y)\}_{y\in \Alive(t)},
	}
where $I(y)$ and $H(y)$ are defined above and $R_t(y)=R_0(y)-t,\, 0\le t\le T_1$, denotes
the remaining lifetime of half-edge $y$. This concludes the initial step in the recursion,
where we defined $\Alive(t)$ and $\SWG(t)$ during the random interval $[T_0,T_1)$.

We continue using induction, by defining  $\Alive(t)$ and $\SWG(t)$ during the random interval $[T_{k},T_{k+1})$, given that the processes are defined on $[0,T_{k})$.
At time $t=T_{k}$, we remove $y^{\star}_{k-1}$ from the set $\Alive(t-)$. 
By construction, we know that $z_k\equiv P_{y^{\star}_{k-1}} \notin \Alive(t-)$, so that $V_{z_k}$ is not a vertex that has been reached by the flow
at time $t$. 
Then, for each of the $d_{V_{z_k}}-1$ other half-edges incident to vertex $V_{z_k}$ we test whether it is \ch{part of a self-loop} or paired to a half-edge from the set 
$\Alive(t-)$. 
\ch{All half-edges incident to $V_{z_k}$  that are part of a self-loop or incident to
$\Alive(t-) $ are removed from vertex $V_{z_k}$; we also remove the involved 
half-edges from the set  $\Alive(t-)$ .}
We will discuss the role of the half-edges incident 
to $V_{z_k}$ that are paired to half-edges in $\Alive(t-)$ in more detail below.

\ch{ For all the remaining siblings of $z_k$ we do the following:
Let $x$ be one such half-edge of $V_{z_k}$}, then $x$ is added
to $\Alive(T_k)$, $I(x)=I(y^{\star}_{k-1})$, $H(x)=H(y^{\star}_{k-1})+1$, while $R_{T_k}(x)$ is an i.i.d.~life time with distribution $G$. We now set $\Alive(t)=\Alive(T_k),\,T_k\le t<T_{k+1}$, where
$T_{k+1}=T_k+\min_{y\in \Alive(T_k)} R_{T_k}(y)$, and where the minimizing half-edge is
called $y^{\star}_k$. Furthermore, for $t\in [T_{k},T_{k+1})$, we can define $\SWG(t)$ by
\eqref{def-flow}, where $R_t(y)=R_{T_k}(y)-(t-T_k)$. Finally, we denote the number 
of \ch{ the $d_{V_{z_k}}-1$ other half-edges incident to vertex $V_{z_k}$ that do not form a self-loop and that are not paired to a half-edge from the set $\Alive(t-)$ by $X_k$.} 
Later, it will also be convenient
to introduce $B_k=d_{V_{z_k}}-1$. Let $S_k=|\Alive(T_k)|$, so that $S_0=X_0^{\sss(1)}+X_0^{\sss(2)}$, 
while $S_k$ satisfies the recursion
    \eqn{
    \label{resurs-S}
    S_k=S_{k-1}+X_k-1.
    }
This describes the evolution of $(\SWG(t))_{t\ge 0}$.

\paragraph{Cycle edges and collision edges.}
At the times \ch{ $T_k,\, k\ge 1,$} we find the half-edge $y_{k-1}^{\star}$ which is paired to $z_k=P_{y_{k-1}^{\star}}$,
and for each of the other half-edges $x$ incident to $V_{z_k}$, we check whether or not $P_x\in \Alive(T_k-)$.
The half-edges paired to alive half-edges in $\Alive(T_k-)$
are special. Indeed, the edge $(x,P_x)$ creates a cycle when $I(x)=I(P_x)$
while $(x,P_x)$ completes a path between $\Ver_1$ and $\Ver_2$,
when $I(x)=3-I(P_x)$. Precisely the latter edges can create the shortest-weight 
path between $\Ver_1, \Ver_2$.  Let us describe these collision 
edges in more detail.

At time $T_k$ and when we create a collision edge consisting of $x_k$ and
$P_{x_k}$, then we record
	\eqn{
	\label{coll-edge}
	\Big((T_k,I(z_k),H(z_k),H(P_{x_k}),R_{T_k}(P_{x_k})\Big)_{k\ge 0}.
	}
It is possible that multiple half-edges incident to $V_{z_k}$ create
collision edges, and if so, we collect all of them in the list in \eqref{coll-edge}.
In this definition it is tempting to write $I(x_k)$ and $H(x_k)$, but note that
$x_k\notin \Alive(T_k)$, whereas its sibbling half-edge $z_k\in \Alive(T_k)$, 
and, moreover, $x_k$ and $z_k$ have the same ancestor and the same height.
\ch{With some abuse of notation we denote  the $i$th collision edge by $(x_i,P_{x_i})$; here $P_{x_i}$ is an alive half-edge and $x_i$ the half-edge which pairs to $P_{x_i}$;
further $z_i$ is the sibling of $x_i$ paired with the minimal edge $y^*$ found by the flow.} 
Let $\Tcol_i$ be the time of creation of the $i$th collision edge.
The weight of the (unique) path between $\Ver_1$ and $\Ver_2$ that passes through the edge
consisting of $x_i$ and $P_{x_i}$ equals $2\Tcol_i+R_{\Tcol_i}(P_{x_i})$, so that the shortest
weight equals:
	\eqn{
	\label{short-weight}
	L_n=\min_{i\ge 0}[2\Tcol_i+R_{\Tcol_i}(P_{x_i})].
	}
Let $I^{\star}$ denote the minimizer of $i\mapsto 2\Tcol_i+R_{\Tcol_i}(P_{x_i})$, then
	\eqn{
	\label{hopcount}
	H_n=H(z_{I^{\star}})+H(P_{x_{I^{\star}}})+1.
	}
Of course, \eqref{short-weight} and \eqref{hopcount} need a proof, which we give now.

{\it Proof that $L_n$ given by \eqref{short-weight} yields the minimal weight.}
Observe that each path between $\Ver_1$ and $\Ver_2$ has a weight $L$ that can be written in the form \ch{$2T_i+R_{T_i}(P_{x_i})$}
 for some $i\ge 0$. Indeed, let $(i_0=\Ver_1,i_1,i_2,\ldots,i_k=\Ver_2)$ form
 a path with weight $L$, and denote the weight on $i_{j-1}i_j$ by $X_{e_j}$ for $1\le j\le k$. For $k=1$, we obviously
find $X_{e_1}=2T_0+X_{e_1}$. For general $k\ge 1$, take the maximal $j\ge 0$ such that
$X_{e_1}+\dots+X_{e_j}\le L/2$. Then, we  write
	$$
	L
	=
	\left\{
	\begin{array}{ll}
	2\sum_{s=1}^j X_{e_s}+[\sum_{s=j+1}^k X_{e_s}-\sum_{s=1}^j X_{e_s}],
	&\text{when }\sum_{s=1}^j X_{e_s}<\sum_{s=j+1}^k X_{e_s},\\
	2\sum_{s=j+1}^k X_{e_s}+[\sum_{s=1}^j X_{e_s}-\sum_{s=j+1}^k X_{e_s}],
	&\text{when }\sum_{s=1}^j X_{e_s}>\sum_{s=j+1}^k X_{e_s},
	\end{array}
	\right.
	$$
which in either case is of the form $L=2T_m+R_{T_m}(y)$, for some  $m\ge 0$ and some
half-edge $y$. Note that in the construction of the flow clusters,
instead of putting weight on the edges, we have given weights to half-edges instead.
In the representation \eqref{short-weight} full edge weight is given to the active half-edges
and weight $0$ to the ones with which they are paired. When the collision edge has been found
we give the full weight to the {\it parent}-edge $P_{x}$. So, in fact, by the redistribution
of the weights in \eqref{short-weight} is an equality in distribution. This completes the proof
of the claim.\hfill \qed

\begin{Remark}[On the number of collision edges]
	\label{rem:finite-pp-enough}
	We do not have to find \emph{all} collision edges. The recursion can be stopped
	when $T_k>L/2$ for some $k\geq 1$, where $L$ is the weight of one the 
	collision edges found previously. This is immediately clear, since all collision 
	edges found at $T_k$ or later have weight exceeding $2T_k>L$.
\end{Remark}

\subsection{Main result: Poisson Point Process limit}
\label{sec-PPP}
In this section, we state our main result, which will imply Theorems \ref{thm-main-first}- \ref{thm-main}. 

\paragraph{Basic constructions and properties.}
To state our main technical result concerning the appearance of collision edges,
we need to define some new constructs. We start by defining a rescaled version of the point
process corresponding to the points in \eqref{coll-edge}. Let us first setup some notation. 
For $i\in \{1,2\}$ and $t\geq 0$, we let
    	\eqn{
    	\label{SWGit-def}
    	|\SWG(t)|=\#\{y\in \Alive(t)\}, \qquad 
	|\SWG^{\sss(i)}(t)|=\#\{y\in \Alive(t)\colon I(y)=i\},
    	}
be the number of alive half-edges at time $t$, as well as those that are closest 
to vertex $i$. By construction, since we check whether the half-edges 
form a cycle or \ch{a} collision edge when the half-edges are born, 
$\SWG^{\sss(1)}(t)$ and $\SWG^{\sss(2)}(t)$ are \emph{disjoint}.
Consider the filtration $(\FF_s)_{s\geq 0}$ with 
$\FF_{s}=\sigma((\SWG(t))_{t\in [0,s]}$ denoting the sigma-algebra generated by
the shortest-weight graph up to time $s$.

Fix a deterministic sequence $s_n\rightarrow \infty$ that will be chosen later on. Now let
    	\eqn{
    	\label{tn-def}
    	t_n=\frac{1}{2\alphan} \log{n},
	\qquad
	\tn=\frac{1}{2\alphan} \log{n}
	-\frac{1}{2\alphan} \log{\big(\WW_{s_n}^{\sss(1)}\WW_{s_n}^{\sss(2)}\big)},
    	}
where, \ch{for $s\geq 0$,}
	\eqn{
	\WW_{s}^{\sss(i)}=\e^{-\alphan s}|\SWG^{\sss(i)}(s)|.
	}
\ch{	
Note that $\e^{\alphan t_n}=\sqrt{n}$, so that at time $t_n$, both $|\SWG^{\sss(i)}(s)|$
are of order $\sqrt{n}$; consequently the variable $t_n$ denotes the typical time at which collision edges start appearing}, and
the time $\tn$ incorporates for stochastic fluctuations in the size of the SWGs. The precise rate at  which $s_n\to\infty$ for asymptotic properties of the construction to hold is determined 
in the proof of Proposition \ref{prop-CTBP-asy} below.
In particular we choose $s_n\rightarrow \infty$ such that $\SWG^{\sss(i)}(t)$ for $t\leq  s_n$ 
can be coupled with two independent two-stage branching processes $\BP^{\sss(i)}(t)$ such that
with high probability $\set{\BP(t) = |\SWG(t)|} ~ \forall t\leq s_n$. 

Define the \emph{residual life-time distribution} $\FR$ to have density $\fR$ given by
	\eqn{
	\label{dens_fR}
	\fR(x)=\frac{\int_0^\infty \e^{-\alpha y}g(x+y)\,dy}{\int_0^\infty \e^{-\alpha y} [1-G(y)]\,dy}.
	}
%

Recall that the $i$th collision edge is given by $(x_i,P_{x_i})$, where $P_{x_i}$ 
is an alive half-edge and $x_i$ the half-edge which pairs to $P_{x_i}$.
In terms of the above definitions, we define
    \eqn{
    \label{resc-variables}
    \barTcol_i=\Tcol_i-\tn,
    \qquad
    \bar{H}_{i}^{\sss({\rm or})}=\frac{H(x_i)- t_n/\bar{\nu}_n}{\sqrt{\bar{\sigma}^2t_n/\bar{\nu}^3}},
    \qquad
    \bar{H}_{i}^{\sss({\rm de})}=\frac{H(P_{x_i})- t_n/\bar{\nu}_n}{\sqrt{\bar{\sigma}^2t_n/\bar{\nu}^3}},
    }
and write the random variables $(\Xi_i)_{i\geq 1}$ with $\Xi_i\in \R\times \{1,2\}\times \R\times \R\times [0,\infty),$
by
    \eqn{
    \label{Ui-def}
    \Xi_i=\big(\barTcol_i, I(x_i), \bar{H}_{i}^{\sss({\rm or})}, \bar{H}_{i}^{\sss({\rm de})}, R_{T_i}(P_{x_i})
    \big).
    }
Then, for sets $A$ in the Borel $\sigma-$algebra of the space $\Scal:= \R\times \{1,2\}\times \R\times \R\times [0,\infty)$,
we define the point process
    \eqn{
    \label{PPP-discr}
    \Pi_n(A)=\sum_{i\geq 1} \delta_{\Xi_i}(A),
    }
where $\delta_x$ gives measure $1$ to the point $x$. Let $\Mcal(\Scal)$ denote the space of all simple locally finite point processes on $\Scal$ equipped with the vague topology (see e.g.\ \cite{KallenBK76}). On this space one can naturally define the notion of  weak convergence of a sequence of random point processes $\Pi_n\in \Mcal(\Scal)$. This is the notion of convergence referred to in the following theorem. In the theorem, we let $\Phi$ denote the distribution function of a standard normal
random variable.

\begin{Theorem}[PPP limit of collision edges]
\label{thm-main-PPP}
Consider the distribution of the point process $\Pi_n \in \Mcal(\Scal)$ 
 defined in
\eqref{PPP-discr} conditional on $(\SWG(s))_{s\in[0,s_n]}$ such that $\WW_{s_n}^{\sss(1)}>0$ and $\WW_{s_n}^{\sss(2)}>0$. 
Then $\Pi_n$  converges in distribution as $n\to\infty$ to a Poisson Point Process (PPP) $\Pi$ 
with intensity measure
    \eqn{
    \label{intensity-PPP}
    \lambda(dt\times i\times dx\times dy\times dr)
    =\frac{2\nu \fR(0)}{\mu}\e^{2\alpha t}dt
    \otimes
    \{1/2,1/2\}
    \otimes
    \Phi(dx)
    \otimes
    \Phi(dy)
    \otimes
    \FR(dr).
    }
\end{Theorem}


\paragraph{Completion of the proof of Theorems \ref{thm-main-first}, \ref{thm-main} and \ref{thm-main-paths}.}
Let us now prove Theorem \ref{thm-main-first} subject to Theorem \ref{thm-main-PPP}.
First of all, by \eqref{resc-variables}, \eqref{short-weight} and \eqref{hopcount} and Remark \ref{rem:finite-pp-enough}, 
    \eqn{
    \Big(\frac{H_n - \frac{1}{\alphan \bar{\nu}}\log{n}}{\sqrt{\frac{\bar{\sigma}^2}{\bar{\nu}^3\alpha}\log{n}}},
    \Wn - \frac{1}{\alphan} \log{n}\Big),
    }
is a continuous function of the point process $\Pi_n$, and, therefore, by the continuous mapping
theorem, the above random variable converges in distribution to some limiting
random variables $(Z,Q)$.

Recall that $I^{\star}$ denotes the minimizer of $i\mapsto 2\Tcol_i+R_{\Tcol_i}(P_{x_i})$.
By \eqref{short-weight}, the weight $\Wn$ as well
as the value of $\Istar$, are functions of the first and the last coordinates of $\Pi_n$.
The hopcount $\Hn$ is a function of the  \ch{third and the fourth}, instead.
By the product form of the intensity in \eqref{intensity-PPP},
we obtain that the limits $(Z,Q)$ are independent.
Therefore, it suffices to study their marginals.
The same observation applies to the multiple path problem in Theorem \ref{thm-main-paths}.

We start with the limiting distribution of the hopcount. By \eqref{resc-variables},
    \eqn{
    \label{lim-decoposition} 
    \frac{H_n - \frac{1}{\alphan \bar{\nu}_n}\log{n}}{\sqrt{\frac{\bar{\sigma}^2}{\bar{\nu}^3\alpha}\log{n}}}
    =\frac{1}{2}\sqrt{2}\bar{H}_{\Istar}^{\sss({\rm or})}+\frac{1}{2}\sqrt{2}\bar{H}_{\Istar}^{\sss({\rm de})}+\op(1).
    }
By Theorem \ref{thm-main-PPP}, the random variables $(\bar{H}_{\Istar}^{\sss({\rm or})},\bar{H}_{\Istar}^{\sss({\rm de})})$, converge to two independent standard normals, so that
also the left-hand side of \eqref{lim-decoposition} 
converges to a standard normal.

The limiting distribution of the weight $\Wn$ is slightly more involved.
By \eqref{tn-def}, 
     \eqref{short-weight} and 
     \eqref{resc-variables},
    	\eqan{
    	\Wn-\frac{1}{\alphan} \log{n}
    	&=\Wn-2t_n=\Wn-2\tn-\frac{1}{\alphan}\log(\WW^{\sss(1)}_{s_n}\WW^{\sss(2)}_{s_n})\\
	&=-\frac{1}{\alphan}\log(\WW^{\sss(1)}_{s_n}\WW^{\sss(2)}_{s_n})
	+\min_{i\geq 1} [2\Tcol_i+R_{\Tcol_i}(P_{x_i})]-2\tn\nn\\
	&=-\frac{1}{\alphan}\log(\WW^{\sss(1)}_{s_n}\WW^{\sss(2)}_{s_n})+\min_{i\geq 1} 
	[2\barTcol_i+R_{\Tcol_i}(P_{x_i})].\nn
    	}
By Proposition \ref{prop-CTBP-asy} below, $(\WW^{\sss(1)}_{s_n},\WW^{\sss(2)}_{s_n})\convd (\WW^{\sss(1)},\WW^{\sss(2)})$,
which are two independent copies of the random variable in \eqref{eqn:tilde-w}. Hence,
    \eqn{
    \Wn-\frac{1}{\alphan} \log{n}\convd
-\frac{1}{\alpha}\log(\WW^{\sss(1)}\WW^{\sss(2)})+
    \min_{i\geq 1} [2P_i+R_i],
    }
where $(P_i)_{i\geq 1}$ form a PPP with intensity $\frac{2\nu\fR(0)}{\mu}\e^{2\alpha t}dt$, and $(R_i)_{i\geq 1}$ are i.i.d.\ random variables with distribution function $\FR$ independently
of $(P_i)_{i\geq 1}$.

We next identify the distribution of $M=\min_{i\geq 1}[2P_i+R_i]$.
First, $(2P_i)_{i\geq 1}$ forms a Poisson process with intensity $\frac{\nu\fR(0)}{\mu}\e^{\alpha t}dt$.
According to  \cite[Example 3.3 on page 137]{resnick}
the point process $(2P_i+R_i)_{i\geq 1}$ is a non-homogeneous Poisson process
with mean-measure the convolution of $\mu(-\infty,x]=\int_{-\infty}^x \frac{\nu\fR(0)}{\mu} \e^{\alpha t}\,dt$ and
$\FR$. Hence $\prob(M\ge x)$ equals the Poisson probability of $0$, where the
parameter of the Poisson distribution is $(\mu*\FR)(x)$, so that
	\eqn{
	\prob(M\ge x)=
 	\exp
	\{-\frac{\nu\fR(0)}{\mu} \e^{\alpha x}
	\int_{0}^{\infty}
 	\FR(z) \e^{-\alpha z}\,dz\}.
	}

Let $\Lambda$ have a Gumbel distribution, i.e., $\prob(\Lambda\leq x)=\e^{-\e^{-x}},\, x\in \mathbb{R}$, then
    \eqn{
    \prob(-a \Lambda+b\geq x)=\e^{-\e^{x/a}\e^{-b/a}}.
    }
From the identity:
$$
\frac{\nu\fR(0)}{\mu} \e^{\alpha x}
	\int_{0}^{\infty}
 	\FR(z) \e^{-\alpha z}\,dz=\e^{x/a}\e^{-b/a},
$$
we conclude that if we take  $a=1/\alpha$ and $b=-\alpha^{-1}\log\Big((\nu \fR(0)/\mu) \int_{0}^{\infty}
 	\FR(z) \e^{-\alpha z}\,dz\Big)$, then
	\eqn{
    	\label{rewrite-Gumbel}
    	\min_{i\geq 1}(2P_i+R_i)\stackrel{d}{=}-\alpha^{-1}\Lambda-\alpha^{-1}\log(\nu \fR(0) B/\mu),
    	}
with $B=\int_{0}^{\infty} \FR(z) \e^{-\alpha z}\,dz$.   In the following lemma, we simplify these
constants:

\begin{Lemma}[The constant]
\label{lem-constants}
The constants $B=\int_{0}^{\infty} \FR(z) \e^{-\alpha z}\,dz$ and $\fR(0)$ are given by
	\eqn{
	B=\bar{\nu}/(\nu-1), 
	\qquad
	\fR(0)=\alpha/(\nu-1). 
	}
Consequently, the constant $c$ in the limit variable \eqref{lim-var}  equals
    	\eqn{
    	\label{factor-sores}
    	c=-\log(\nu \fR(0)B/\mu)=
    	\log(\mu(\nu-1)^2/(\alpha\nu\bar{\nu})).
    	}
\end{Lemma}

\proof We start by computing $\fR(0)$, for which we note that by \eqref{dens_fR} and \eqref{eqn:malthusian-param},
	\eqn{
	\fR(0)=\frac{\int_0^\infty \e^{-\alpha y}g(y)\,dy}{\int_0^\infty \e^{-\alpha y} [1-G(y)]\,dy}
	=\Big(\nu \int_0^\infty \e^{-\alpha y} [1-G(y)]\,dy\Big)^{-1}.
	}
Further, by partial integration,
	\eqn{
	\label{int-2-comp}
	\int_0^\infty \e^{-\alpha y} [1-G(y)]\,dy
	=\Big[-\frac{1}{\alpha}\e^{-\alpha y} [1-G(y)]\Big]_{y=0}^{\infty}
	-\frac{1}{\alpha}\int_0^\infty \e^{-\alpha y} g(y)\,dy
	=\frac{1}{\alpha}-\frac{1}{\alpha\nu}=\frac{\nu-1}{\alpha \nu},
	}
where we again use \eqref{eqn:malthusian-param}. Combining both equalities yields $\fR(0)=\alpha/(\nu-1)$.

\ch{For $B$, we again use partial integration, followed by the substitution of \eqref{dens_fR}; this yields
	\eqan{
	B&=\int_{0}^{\infty} \FR(z) \e^{-\alpha z}\,dz=\frac{1}{\alpha}\int_{0}^{\infty} \fR(z) \e^{-\alpha z}\,dz\\
	&=\frac{\nu}{\nu-1}\int_{0}^{\infty}\e^{-\alpha z} \int_0^{\infty} \e^{-\alpha y} g(y+z)\,dy \,dz,\nn
	}
}
by \eqref{int-2-comp}. The final integral can be computed using\ch{
\begin{eqnarray}
&&	\int_{-\infty}^\infty\e^{-\alpha z}\indic{z\ge 0}\int_{-\infty}^{\infty} \e^{-\alpha y} g(y+z)\indic{y\ge 0}\,dy \,dz\nonumber\\
&&\qquad	=\int_0^{\infty} sg(s)\e^{-\alpha s}\,ds=\frac{1}{\nu} \int_0^{\infty} s\bar{G}(ds)=\bar{\nu}/\nu.
\end{eqnarray}
}
\ch{This completes the proof of both Theorem \ref{thm-main-first} and Theorem \ref{thm-main}, given Theorem \ref{thm-main-PPP}}.
\hfill\qed

\subsection{Overview of the proof of Theorem \ref{thm-main-PPP}}
\label{sec-overview}
In this section, we reduce the proof of Theorem \ref{thm-main-PPP} to two 
main propositions. Recall the shortest weight graph or flow cluster
$\SWG(t)$ defined in the previous section as well as the associated filtration $(\FF_t)_{t\geq 0}$.
We shall couple these flow clusters from two points with $(\BP(t))_{t\geq 0}$ 
where $\BP(t)=(\BP^{\sss(1)}(t),\BP^{\sss(2)}(t))$ are two independent CTBPs
starting with offspring distribution $D$.  For a prescribed such coupling of $(\SWG(t))_{t\geq 0}$ and
$(\BP(t))_{t\geq 0}$, we let $\SWG(t)\symdiff\BP(t)$ denote the set 
\ch{of} miscoupled half-edges \ch{at time $t$}. Then we prove the following limiting result:

\begin{Proposition}[Coupling the SWG to a BP]
\label{prop-CTBP-asy}
(a) There exists $s_n\rightarrow \infty$ and a coupling of $(\SWG(s))_{s\geq 0}$ and $(\BP(s))_{s\geq 0}$
such that
	\eqn{
	\prob\Big((\SWG(s))_{s\in[0,s_n]}= (\BP(s))_{s\in[0,s_n]}
	\Big)=1-o(1).
	}
Consequently, with $\WW_{s_n}^{\sss(i)}=\e^{-\alphan s_n}|\SWG^{\sss(i)}(s_n)|$,
	\eqn{
	\label{WW-vep}
	\liminf_{\vep\downarrow 0}
	\liminf_{n\ra \infty}
	\prob\Big( \WW_{s_n}^{\sss(1)}\in [\vep, 1/\vep],
	\WW_{s_n}^{\sss(2)}\in [\vep, 1/\vep]
	\Big| \WW_{s_n}^{\sss(1)}>0, \WW_{s_n}^{\sss(2)}>0\Big)=1.
	}
(b) There exists a coupling of  $(\SWG(s))_{s\geq 0}$ and 
$(\BP_{\sss(n)}(s))_{s\geq 0}$, 
and sequences $\vep_n=o(1)$ and $\Bn\ra \infty$  
such that, conditionally on $\FF_{s_n}$,  
    	\eqn{
    	\prob\Big(|\SWG(t_n+\Bn)\symdiff\BP_{\sss(n)}(t_n+\Bn)| \geq \vep_n \sqrt{n}\mid \FF_{s_n}\Big)
    	\convp 0,
    	}
where $(\BP_{\sss(n)}(t)_{t\geq 0}=(\BP^{\sss(1)}_{\sss(n)}(t),\BP^{\sss(2)}_{\sss(n)}(t))_{t\geq 0}$
and 
$(\BP^{\sss(1)}_{\sss(n)}(s))_{s\geq s_n}, (\BP^{\sss(2)}_{\sss(n)}(s))_{s\geq s_n}$ are two independent
two-stage Bellman-Harris processes with 
offspring \ch{ $D^{\star}_n-1$ (where $D^{\star}_n$ has the sixe-biased distribution $F^{\star}_n$ of
$F_n$, see \eqref{eqn:size-biased})} for every individual, 
and edge weights with continuous distribution function $G$, 
and starting at time $s_n$ in $\BP(s_n)$ from
part (a), respectively. 
\end{Proposition}

The proof of Proposition \ref{prop-CTBP-asy} is deferred to Section \ref{sec-thinning-CM}.
\ch{In the sequel, we shall assume that $\prob$ denotes the coupling measure from Proposition \ref{prop-CTBP-asy}. In particular, this yields a coupling between $\CMnd$ for different $n\geq 1$,
as well as a coupling between $\CMnd$ and the $n$-dependent branching processes 
$(\BP(s))_{s\geq 0}$. Under this coupling law, we can speak of convergence in probability,
and we shall frequently do this in the sequel.}

For $i\in \{1,2\}$, $k\geq 0$, and $t\geq 0$, we define
    \eqn{
    \label{SWGit[t,t+s)-def}
    |\SWG^{\sss(i)}_k[t,t+s)|=\#\{y\in \Alive(t)\colon I(y)=i, H(y)=k, R_t(y)\in [0,s)\},
    }
as the number of alive half-edges at time $t$ that (a) are in the SWG of vertex $\Ver_i$,
(b) have height $k$, and (c) have remaining lifetime at most $s$. We
further write
	\eqn{
	\label{SWGleqkt[t,t+s)-def}
    	|\SWG^{\sss(i)}_{\sss \leq k}[t,t+s)|=\#\{y\in \Alive(t)\colon I(y)=i, H(y)\leq k, R_t(y)\in [0,s)\},
    	}
for the number of vertices that have height at most $k$.
To formulate the CLT for the height of vertices, we will choose
    \eqn{
    \label{ktx-def}
    k_t(x) = \frac{t}{\bar{\nu}} + x \sqrt{t\frac{\bar{\sigma}^2}{\bar{\nu}^3}}.
    }
Finally, for a half-edge $y\in \Alive(t)$, we let $X_y^{\star}=d_{V_y}-1$. 

\begin{Proposition}[Ages and heights in SWG]
\label{prop-CLT-stable-age}
Fix $j\in \{1,2\}$, $x,y,t \in \Rbold$, $ s_1,s_2>0$. Then, conditionally on
$\WW_{s_n}^{\sss(1)}\WW_{s_n}^{\sss(2)}>0$,\\
(a)
    	\eqan{
    	&\e^{-2\alphan t_n}
	|\SWG_{\sss \leq k_{t_n}(x)}^{\sss(j)}[\tn+t,\tn+t+s_1)|
	|\SWG_{\sss \leq k_{t_n}(y)}^{\sss(3-j)}[\tn+t,\tn+t+s_2)|\\
	&\qquad\qquad\convp \e^{2\alpha t}\Phi(x)\Phi(y)\FR(s_1)\FR(s_2),\nn
    	\label{eqn:ratio-k-normal}
    	}
(b)
	\eqan{
	&\e^{-2\alphan t_n}
	|\SWG_{\sss \leq k_{t_n}(x)}^{\sss(j)}[\tn+t,\tn+t+s_1)|
	\sum_{v} X_v^{\star}\1_{\{v\in \SWG_{\sss \leq k_{t_n}(y)}^{\sss(3-j)}[\tn+t,\tn+t+s_2)\}}\\
	&\qquad\qquad\convp \nu\e^{2\alpha t}\Phi(x)\Phi(y)\FR(s_1)\FR(s_2).\nn
    \label{eqn:ratio-k-normal-sec}
    }
\end{Proposition}
The first assertion in the above proposition follows from
\cite[Theorem 1(b)]{samuels1971distribution} in the case that our
CTBP has finite-variance offspring. The proof of Proposition \ref{prop-CLT-stable-age}
is deferred to Section \ref{sec-CTBP-pf-strong}.

\paragraph{Completion of the proof of Theorem \ref{thm-main-PPP}.}
Recall that $\FF_t=\sigma((\SWG(s))_{s\in[0,t]})$. We will investigate the number of collision
edges $(x_i,P_{x_i})$ with $I(x_i)=j\in \{1,2\}$, $H(x_i)\leq k_{t_n}(x)$, $H(P_{x_i})\leq k_{t_n}(y)$ and
$R_{\Tcol_i}(P_{x_i})\in [0,s)$ created in the time interval $[\tn+t,\tn+t+\vep)$, where $\vep>0$ is  small.
We let $\II=[a,b)\times \{j\}\times (-\infty,x]\times (-\infty, y]\times [0,s]$ be a subset of $\Scal$, 
and we prove that
	\eqn{
	\label{aim-PPP}
	\prob(\Pi_n(\II)=0\mid \FF_{s_n})\convp 
	\exp{\big\{-\int_a^b \frac{2\nu \fR(0)}{\mu}\e^{2\alpha t}\Phi(x)\Phi(y)\FR(s)dt\big\}}.
	}
By \cite[Theorem 4.7]{KallenBK76}, this proves the claim.

We split
	\eqn{
	\II=\bigcup_{\ell=1}^N \II_\ell^{\sss(\vep)},
	}
where $\II_l^{\sss(\vep)}=[t_{\ell-1}^{\sss(\vep)}, t_{\ell}^{\sss(\vep)})\times \{j\}\times (-\infty,x]\times (-\infty, y]\times [0,s),$ 
with $t_{\ell}^{\sss(\vep)}=a+\ell\vep$ and $\vep=(b-a)/N$, with $N\in \Nbold$. 
We will let $\vep\downarrow 0$ later on.
For a fixed $\vep>0$, we say that a collision edge $(x_i,P_{x_i})$ is a \emph{first round collision edge}
when there exists $j\in[N]$ and a half-edge $y\in \Alive(t_{l-1}^{\sss(\vep)})$ such that
$y$ is found by the flow in the time interval $\II_\ell^{\sss(\vep)}$, $y$ is paired to the half-edge 
$P_y$ whose sibling half-edge $x_i$ is paired to $P_{x_i}\in \Alive(t_{\ell-1}^{\sss(\vep)})$
with $I(y)=j\neq I(P_{x_i})=3-j$. We call all other collision edges \emph{second round collision edges.}
Denote the point processes of first and second round collision edges by $\PiFR_n$ and $\PiSR_n$,
so that $\Pi_n=\PiFR_n+\PiSR_n$. The next two lemmas investigate the point processes
$\PiFR_n$ and $\PiSR_n$:

\begin{lemma}[PPP limit for the first round collision edges]
\label{lem-FR-coll}
For every $s\geq 0$, $x,y\in \R$, $j\in \{1,2\}$, $\vep>0$ and $\ell\in[N]$, as $n\ra \infty$,
	\eqn{
	\prob\big(\PiFR_n(\II_\ell^{\sss(\vep)})=0\mid \FF_{t_{\ell-1}^{\sss(\vep)}}\big)\convp 
	\exp{\big\{-\e^{2\alpha t_{\ell-1}^{\sss(\vep)}}\Phi(x)\Phi(y)\FR(s)\FR(\vep)\big\}}.
	}
\end{lemma}

\proof The number of half-edges $z\in \Alive(\tn+t_{\ell-1}^{\sss(\vep)})$ that are found by the flow having $I(z)=j$ and $H(z)\leq k_{t_n}(x)$
is equal to
    \eqn{
    \label{step-1}
    |\SWG_{\sss \leq k_{t_n}(x)}^{\sss(j)}[\tn+t_{\ell-1}^{\sss(\vep)},\tn+t_{\ell-1}^{\sss(\vep)}+\vep)|.
    }
Fix such a half-edge $z$, and note that it is paired to $P_z$ that has $X^{\star}_z=d_{V_{P_z}}-1$
sibling half-edges. For each of these half-edges we test whether it is paired to a half-edge
in $\Alive(\tn+t_{\ell-1}^{\sss(\vep)})$ or not. Therefore, the total number of tests performed in
the time interval $[t_{\ell-1}^{\sss(\vep)}, t_{\ell}^{\sss(\vep)})$ is equal to
    \eqn{
    \label{step-2}
    \sum_{z}X_z^{\star}\1_{\{z\in \SWG_{\sss  \leq k_{t_n}(x)}^{\sss(j)}[\tn+t_{\ell-1}^{\sss(\vep)},\tn+t_{\ell}^{\sss(\vep)})\}}.
    }
By construction, we test whether these half-edges are paired to half-edges that are
incident to the SWG or not. Each of these edges is paired to a 
half-edge $w\in \Alive(\tn+t_{\ell-1}^{\sss(\vep)})$ with $I(w)=3-j$
(and thus creating a collision edge) and \ch{$H(w)\leq k_{t_n}(y)$} and
\ch{$R_{\tn+t_{\ell-1}^{\sss(\vep)}}(w)\in [0,s)$} with probability equal to
\ch{    \eqn{
    \label{step-3}
    \frac{1}{\ell_n-o(n)}  |\SWG_{\sss \leq k_{t_n}(y)}^{\sss(3-j)}[\tn+t_{\ell-1}^{\sss(\vep)},\tn+t_{\ell-1}^{\sss(\vep)}+s)|.
    }
    }
Therefore, the expected number of first round collision edges
$(x_i,P_{x_i})$ with $I(x_i)=j\in \{1,2\}$, $H(x_i)\leq k_{t_n}(x)$, $H(P_{x_i})\leq k_{t_n}(y)$ and
$R_{\tn+t_{\ell-1}^{\sss(\vep)}}(P_{x_i})\in [0,s)$ created in the time interval $[\tn+t_{\ell-1}^{\sss(\vep)}, \tn+t_{\ell}^{\sss(\vep)})$ 
equals  \ch{the product of the expressions in \eqref{step-2} and \eqref{step-3}, and can be 
rewritten as}
    \ch{
    \eqan{
    \label{step-4}
      &\frac{1}{\ell_n-o(n)} \e^{2\alphan t_n}
    \Big(\e^{-\alphan t_n}\sum_z X_z^{\star}
	\1_{\{z\in \SWG_{\sss \leq k_{t_n}(x)}^{\sss(j)}[\tn+t_{\ell-1}^{\sss(\vep)},\tn+t_{\ell}^{\sss(\vep)})\}}\Big)\nn\\
	&\qquad \times
    \Big(\e^{-\alphan t_n}|\SWG_{\sss \leq k_{t_n}(y)}^{\sss(3-j)}[\tn+t_{\ell-1}^{\sss(\vep)},\tn+t_{\ell-1}^{\sss(\vep)}+s)|\Big).
    }
    }
\ch{
By Proposition \ref{prop-CLT-stable-age}, conditionally on $\FF_{\tn+t_{\ell-1}^{\sss(\vep)}}$,
and using that
$
(\ell_n-o(n))^{-1} \e^{2\alphan t_n}\to \mu^{-1},
$
}
\ch{ we find that \eqref{step-4}, which represents the expected number of collision edges
$x_i$ with $I(x_i)=j\in \{1,2\}$, $H(x_i)\leq k_{t_n}(x)$, $H(P_{x_i})\leq k_{t_n}(y)$ and
$R_{\tn+t_{\ell-1}^{\sss(\vep)}}(P_{x_i})\in [0,s)$ created in the time interval $[\tn+t_{\ell-1}^{\sss(\vep)}, \tn+t_{\ell}^{\sss(\vep)})$, converges in probability to:}
    \eqn{
    \label{step-5}
    \frac{\nu}{\mu}\e^{2\alpha t}\Phi(x)\Phi(y)\FR(s) \FR(\vep).
    }
Further, for $\vep>0$, conditionally on $\FF_{\tn+t_{\ell-1}^{\sss(\vep)}}$, the probability that
none of the half-edges found in the time interval in between 
$[\tn+t_{\ell-1}^{\sss(\vep)},\tn+t_{\ell}^{\sss(\vep)})$ creates a collision edge is asymptotically equal to
	\eqan{
	&
	\prod_{v\in \SWG_{\leq k_{\tn}(x)}^{\sss(j)}[\tn+t_{\ell-1}^{\sss(\vep)},\tn+t_{\ell}^{\sss(\vep)})}
	\Big(1-\frac{1}{\ell_n-o(n)}
	|\SWG_{\sss \leq k_{t_n}(y)}^{\sss(3-j)}[\tn+t_{\ell-1}^{\sss(\vep)},\tn+t_{\ell-1}^{\sss(\vep)}+s)|\Big)^{X_v^{\star}}\\
	&\qquad\quad\qquad\convp \exp\big\{-\frac{\nu}{\mu}\e^{2\alpha t_{\ell-1}^{\sss(\vep)}}\Phi(x)\Phi(y)\FR(s) \FR(\vep)\big\}.\nn
	}
\hfill\qed

\begin{lemma}[A bound on the second round collision edges]
\label{lem-SR-coll}
\ch{For  $x,y\in \R$, $j\in \{1,2\}$, $\vep>0$ and $\ell\in[N]$, as $n\ra \infty$,}
	\eqn{
	\prob\big(\PiSR_n(\II_\ell^{\sss(\vep)})\geq 1\mid \FF_{t_{\ell-1}^{\sss(\vep)}}\big)=\Op(1)\FR(\vep)G(\vep).
	}
\end{lemma}

\proof
By analogous arguments as above, the expected number of second round collision edges
is of order
	\eqn{
	\Op(1)\e^{2\alpha t}\Phi(x)\Phi(y)\FR(s) \FR(\vep) G(\vep),
    	}
since one of the half-edges $z$ that is found by the flow in the time interval
$[\tn+t_{\ell-1}^{\sss(\vep)},\tn+t_{\ell}^{\sss(\vep)})$ needs to satisfy that one of the $d_{V_{P_z}}-1$ 
half-edges has weight at most $\vep$, and which, upon being found,
needs to create a collision edge. 
\hfill\qed
\medskip

Now we are ready to complete the proof of Theorem \ref{thm-main-PPP}. We use that
	\eqan{
	\prob\big(\Pi_n(\II)=0\mid \FF_{s_n}\big)
	=\expec\Big[\prod_{\ell=1}^{N} \prob\big(\Pi_n(\II_\ell^{\sss(\vep)})=0\mid \FF_{t_{\ell-1}^{\sss(\vep)}}\big)\mid \FF_{s_n}\Big].
	}
We start with the upper bound, for which we use that
	\eqan{
	\prob\big(\Pi_n(\II_\ell^{\sss(\vep)})=0\mid \FF_{t_{\ell-1}^{\sss(\vep)}}\big)
	&\leq \prob\big(\PiFR_n(\II_\ell^{\sss(\vep)})=0\mid \FF_{t_{\ell-1}^{\sss(\vep)}}\big)\\
	&\convp \exp{\big\{-\e^{2\alpha t_{\ell-1}^{\sss(\vep)}}\Phi(x)\Phi(y)\FR(s)\FR(\vep)\big\}},\nn
	}
by Lemma \ref{lem-FR-coll}. We conclude that
	\eqan{
	\prob\big(\Pi_n(\II)=0\mid \FF_{s_n}\big)
	&\leq \expec\Big[\prod_{\ell=1}^{N} 
	\exp{\big\{-\e^{2\alpha t_{\ell-1}^{\sss(\vep)}}\Phi(x)\Phi(y)\FR(s)\FR(\vep)\big\}}\mid \FF_{s_n}\Big]\\
	&=
	\exp{\big\{-\sum_{\ell=1}^N \e^{2\alpha t_{\ell-1}^{\sss(\vep)}}\Phi(x)\Phi(y)\FR(s)\FR(\vep)\big\}}\nn\\
	&\ra \exp{\big\{-\fR(0)\int_a^b \e^{2\alpha t}\Phi(x)\Phi(y)\FR(s)dt\big\}},\nn
	}
since $\lim_{\vep\downarrow 0} \FR(\vep)/\vep =\fR(0)$, and the Riemann approximation
	\eqn{
	\vep\sum_{\ell=1}^N \e^{2\alpha t_{\ell-1}^{\sss(\vep)}}\ra \int_a^b \e^{2\alpha t}dt.
	}
This proves the upper bound. 

For the lower bound, we instead bound
	\eqan{
	\prob\big(\Pi_n(\II)=0\mid \FF_{s_n}\big)
	&\geq \expec\Big[\prod_{\ell=1}^{N} \prob\big(\PiFR_n(\II_\ell^{\sss(\vep)})=0\mid \FF_{t_{\ell-1}^{\sss(\vep)}}\big)\mid \FF_{s_n}\Big]\\
	&\qquad 
	-\expec\Big[\Big(\sum_{\ell=1}^{N} \prob\big(\PiSR_n(\II_\ell^{\sss(\vep)})\geq 1
	\mid \FF_{t_{\ell-1}^{\sss(\vep)}}\big)\Big) \wedge 1\mid \FF_{s_n}\Big].\nn
	}
The first term has already been dealt with, the second term is, by Lemma \ref{lem-SR-coll}, bounded by
	\eqn{
	\expec\Big[\Big(\Op(1)\sum_{\ell=1}^{N}\FR(\vep)G(\vep)\Big)\wedge 1\mid \FF_{s_n}\Big]
	=\op(1) ,
	}
as $\vep\downarrow 0$, by dominated convergence, since $\FR(\vep)=\vep \fR(0)(1+o(1))$ and $G(\vep)=o(1)$.
\hfill\qed
%

\section{Height CLT and stable age for CTBP}
\label{sec-CTBP-pf-strong}

In this section, we set the stage for the proof of Proposition \ref{prop-CLT-stable-age}
for CTBPs, by investigating the first and second moment of particles of several types.
We will make use of second moment methods similar to the ones in
\cite{samuels1971distribution}, but with a suitable truncation argument
to circumvent the problem of infinite-variance offspring distributions.

We take $K\geq 1$ large and define, for an appropriate $\eta\in (0,1)$
that will be determined later on,
    \eqn{
    \label{miK-def}
    m_i=K\eta^{-i},
    }
and investigate the Bellman-Harris process where each individual in generation
$i$ has offspring distribution $(X\wedge m_i)$ instead of $X$,
\ch{where $X$ denotes the offspring of our CTBP}.

We denote the number of alive individuals in generation $k$ at time $t$
in the original branching process by $|\BP_k(t)|$, and let $|\BP_k[t,t+s)|$
denote the number of alive individuals in generation $k$ with residual lifetime
at most $s$. We let $|\BP_k^{\sss (\vec{m})}[t,t+s)|$
denote the number of individuals in generation $k$ at time $t$ and with remaining lifetime
at most $s$ of the \emph{truncated} branching
process. Define
    \eqn{
    |\BP_{\sss \leq k}^{\sss (\vec{m})}[t,t+s)|=\sum_{j=0}^k |\BP_j^{\sss (\vec{m})}[t,t+s)|,
    \qquad
    |\BP^{\sss (\vec{m})}[t,t+s)|=\sum_{j=0}^\infty |\BP_j^{\sss (\vec{m})}[t,t+s)|.
    }
We also write $|\BP^{\sss (\vec{m})}(t)|=\lim_{s\rightarrow \infty} |\BP^{\sss (\vec{m})}[t,t+s)|$. 
A key ingredient to the proof of Proposition \ref{prop-CLT-stable-age} is 
Proposition \ref{prop-sec-mom-clt} below. In its statement, we also use $\eta=\nu \int_0^{\infty} \e^{-2\alpha s}\,d{\GF}(s)$,
so that $\eta<1$ since $\alpha$ is such that $\nu \int_0^{\infty} \e^{-\alpha s}\,d{\GF}(s)=1$.

\begin{Proposition}[First and second moment CLT]
\label{prop-sec-mom-clt}
\ch{Choose $m_i=K\eta^{-i}$ as in \eqref{miK-def}.} Assume that
the $X\log{X}$ condition holds, i.e., $\expec[X\log(X)_+]<\infty$ where $X$ is the random amount of offspring of our CTBP. Then with $A=(\nu-1)/\alpha\nu{\bar \nu}$,\\
(a) 
    	\eqn{
    	\lim_{t\ra \infty}\e^{-\alpha t}\expec\big[|\BP(t)|\big]=A,
	\qquad \lim_{t\ra \infty}\e^{-\alpha t}\expec\Big[|\BP(t)|-|\BP^{\sss (\vec{m})}(t)|\Big]
    	=0,
    	}
(b) there exists a $C>0$ such that uniformly in $t\rightarrow \infty,$
    \ch{	\eqn{
    	\e^{-2\alpha t}\expec\big[|\BP^{\sss (\vec{m})}(t)|^2\big]
	\leq C K,
	}
	}
(c)
    	\eqn{
	\label{lim-trun-case}
    	\lim_{t\ra \infty}\e^{-\alpha t}\expec\big[|\BP_{\sss \leq k_t(x)}^{\sss (\vec{m})}[t,t+s)|\big]
	=A\Phi(x)\FR(s),
    	}
where $k_t(x)$ is defined in \eqref{ktx-def}.\\
(d) The same results hold uniformly in $n$ when $t=t_n=\frac{1}{2\alphan}\log{n}$ as in \eqref{tn-def},
$\alpha$ is replaced by $\alphan$, $K$ by $K_n$
and the branching process offspring distribution $X_n$ depends on $n$ in such a way that
$X_n\convd X$, $\expec[X_n]\to \expec[X]$ and $\limsup_n \expec[X_n \log(X_n/\Kn)_+]=0$,
for any $\Kn\rightarrow \infty$.
\end{Proposition}

\proof
We start by proving Proposition \ref{prop-sec-mom-clt}(a). The first claim is proved in \cite{jagers1975branching, jagers-nerman4}. 
We bound the first moment of the difference between the
truncated and the original branching process. Let $\nu$ be the expected offspring of the
Bellman-Harris process, and let $\nu^{\sss(i)}=\expec[X\indic{X\leq m_i}]$, where $m_i=K\eta^{-i}$. We compute that
    	\eqan{
	\label{exp-Zt}
    	\e^{-\alpha t}\expec\Big[\sum_{k=0}^{\infty} [\BP_k(t)|-|\BP_k^{\sss (\vec{m})}(t)|]\Big]
	&=\e^{-\alpha t} \sum_{k=0}^{\infty} [\nu^k-\prod_{i=1}^k \nu^{\sss(i)}] \big[G^{\star k}(t)-G^{\star (k+1)}(t)\big],
    	}
where $G$ is the distribution function of the edge weights.
In order to bound the differences $\nu^k-\prod_{i=1}^k \nu^{\sss(i)}$, we rely on the following lemma:

\begin{Lemma}[Effect of truncation on expectation CTBP]
\label{lem-trunc-eff}
Fix $\eta\in (0,1)$ and $m_i=K\eta^{-i}$, if $\expec[X\log{(X)}_+]<\infty$, then
    \eqn{
    [1-\prod_{i=1}^\infty \frac{\nu^{\sss(i)}}{\nu}]\leq (\log(1/\eta))^{-1}
	\expec\big[X\log{(X/K)}_+\big]=o_{\sss K}(1),
    }
where $o_{\sss K}(1)$ denotes a quantity that converges to zero as
$K\rightarrow \infty$.
\end{Lemma}

\proof  Since $\frac{\nu^{\sss(i)}}{\nu}\le 1$ for all $i\geq 1$, it is easily shown by induction that
    	\eqn{
    	1-\prod_{i=1}^k \frac{\nu^{\sss(i)}}{\nu}
	\leq
    	\sum_{i=1}^{k} (1-\frac{\nu^{\sss(i)}}{\nu})
    \leq
    \sum_{i=1}^{\infty} (1-\frac{\nu^{\sss(i)}}{\nu}).
}
Now, using that $\nu>1$,
    	\eqn{
    	\sum_{i=1}^\infty\big(1-\frac{\nu^{\sss(i)}}{\nu}\big)
    	\leq\sum_{i=1}^{\infty} \expec[X\indic{X>m_i}]
    	=\expec\big[X\sum_{i=1}^{\infty} \indic{m_i<X}\big],
    	}
and we note that the number of $i$ for which $m_i=K\eta^{-i}\leq x$ is
at most $[\log{(x/K)}/\log{(1/\eta)}]\vee 0$.
Therefore
    	\eqn{
	\label{prod-sum-ineq}
    	1-\prod_{i=1}^k \frac{\nu^{\sss(i)}}{\nu}
	\leq \sum_{i=1}^{\infty} (1-\frac{\nu^{\sss(i)}}{\nu})
	\leq (\log(1/\eta))^{-1}
	\expec\big[X\log{(X/K)}_+\big],
    	}
which converges to zero when $K\rightarrow \infty$.
\hfill\qed
\medskip

By Lemma \ref{lem-trunc-eff} and \eqref{exp-Zt},
    \eqn{
    \label{smalloK}
    \e^{-\alpha t}\expec\Big[\sum_{k=0}^{\infty} [|\BP_k(t)|-|\BP_k^{\sss (\vec{m})}(t)|]\Big]
    =o_{\sss K}(1)\e^{-\alpha t}\expec\Big[\sum_{k=0}^{\infty} |\BP_k(t)|\Big]
    =o_{\sss K}(1),
    }
which completes the proof of Proposition \ref{prop-sec-mom-clt}(a).


We continue with the proof of the second moment estimate in Proposition \ref{prop-sec-mom-clt}(b).
We follow the proof in \cite{samuels1971distribution},
keeping attention to the truncation.  We introduce $h$ as the generating function of $X$, and $h_j$ as the generating function
of $(X\wedge m_j)$, i.e.,
	\eqn{
	h(s)=\expec[s^{X}],\qquad h_j(s)=\expec[s^{(X\wedge m_j)}],
	}
where $m_j$ is given by \eqref{miK-def}. Parallel to calculations in the proof of \cite[Lemma 4]{samuels1971distribution},
	\begin{eqnarray}
	\label{recursiontweede1}
	&&\expec[|\BP^{\sss (\vec{m})}|^2]=h_1''(1)(\expec[|\BP^{\sss (\vec{m_1})}|])^2*\GF
	+h_1'(1)\expec[|\BP^{\sss (\vec{m_1})}|^2]*\GF,
	\end{eqnarray}
\ch{where $\vec{m}_1=(m_2, m_3, \ldots)$, is $\vec{m}$ with the first element removed.}
Transforming to
	\eqn{
	|\BPbar^{\sss (\vec{m})}(t)|=\e^{-\alpha t}|\BP^{\sss (\vec{m})}(t)|,
	}
we obtain, by multiplying both sides of \eqref{recursiontweede1}  by $\e^{-2\alpha t}$,
	\begin{eqnarray}
	\label{rectrans1}
	&&\expec[|\BPbar^{\sss (\vec{m})}|^2]
	=\frac{\eta h_1''(1)}{\nu}(\expec[|\BPbar^{\sss (\vec{m_1})}|])^2*Q+
	\frac{\eta h_1'(1)}{\nu}\expec[|\BPbar^{\sss (\vec{m_1})}|^2]*{\bar \GF},
	\end{eqnarray}
where\ch{
	\eqn{
	\label{tranforms}
	{\bar \GF}(t)=\nu\int_0^t \e^{-\alpha y}\,d\GF(y),
	\qquad
	Q(t)=\eta^{-1}\int_0^t \e^{-\alpha y}\,d{\bar \GF}(y)=
	\eta^{-1}\nu \int_0^t \e^{-2\alpha y}\,d\GF(y),
	}
	}
and where we recall that $\eta=\int_0^{\infty} \e^{-\alpha y}\,d{\bar \GF}(y)<1$ and $\nu=h'(1)$.
Iteration of \eqref{rectrans1} yields
	\eqan{
	\label{driezeventien}
	\expec[|\BPbar^{\sss (\vec{m})}|^2]
	&=\sum_{j=1}^\infty b_1\cdots b_{j-1}a_j \expec[|\BPbar^{\sss (\vec{m}_j)}|]^2*{\bar \GF}_{j-1}*Q,
	}
where
	\eqn{
	a_j=\frac{\eta h''_j(1)}{\nu},\qquad
	b_j=\frac{\eta h'_j(1)}{\nu},
	}
and \ch{where $\vec{m}_j=(m_{j+1}, m_{j+2}, \ldots)$.}	
According to \eqref{smalloK} the expectation $\expec[|\BPbar^{\sss (\vec{m}_j)}(t_n)|]$ has the same asymptotic behavior as $\expec[\BPbar(t_n)|]$. Hence, by \cite[Lemma 1(a)]{samuels1971distribution} or alternatively by  part (a),
	\eqn{
\label{iter-second-moment}
	\lim_{n\to \infty} \expec[|\BPbar^{\sss (\vec{m}_j)}(t_n)|]
=\lim_{n\to \infty} \expec[|\BPbar(t_n)|]= A.
	}
Since $b_1b_2\ldots b_j\leq \eta^j$ tends to zero exponentially, this leads to
	\eqn{
	\label{E[Z2]}
	\lim_{n\to \infty}
	\expec[|\BPbar^{\sss (\vec{m})}(t_n)|^2]
	=
	A^2
	\sum_{j=1}^\infty b_1\cdots b_{j-1}a_j.
	}
We bound the arising sum in the following lemma:

\begin{Lemma}[Effect of truncation on variance CTBP]
\label{lem-trunc-eff-var}
For $m_i=K\eta^{-i}$, and with $\nu=\expec[X]$,
    \eqn{
    \sum_{j=1}^\infty b_1\cdots b_{j-1}a_j\leq \frac{2\nu K}{1-\eta}.
    }
\end{Lemma}


\proof We bound $ b_j\leq \eta$, and 

    \eqn{
    \label{sigma-bnd}
    a_{j}
    \leq \eta \expec[(X\wedge m_j)^2]
    =\eta\Big(m_j^2 \prob(X>m_j)+\expec[X^2 \indic{X\leq m_j}]\Big),
    }
so that
  	\eqn{
  	\sum_{j=1}^\infty b_1\ldots b_{j-1}a_j\leq
  	\sum_{j=1}^\infty m_j^2 \prob(X>m_j)\eta^j+
  	\sum_{j=1}^\infty  \expec[X^2 \indic{X\leq m_j}]\eta^j.
  }
We bound both terms separately. The first contribution equals
    \eqn{
    K^2\sum_{j=1}^{\infty}\prob(X>K\eta^{-j}) \eta^{-j}
    =K^2 \expec[\sum_{j=1}^{\infty}\eta^{-j}\indic{K\eta^{-j}<X}]
    =K^2\expec[\frac{\eta^{-a(X)}-1}{1-\eta}],
    }
where $a(x)=\max\{j\colon K\eta^{-j}<x\}=\lfloor \log{(x/K)}/\log{(1/\eta)}\rfloor$. Therefore, $\eta^{-a(X)}\leq X/K$,
so that
    \eqn{
    \sum_{j=1}^{\infty}m_j^2 \prob(X>m_j) \eta^{j}
    \leq \frac{K^2}{1-\eta}\expec[X/K]=\frac{K\nu}{1-\eta}.
    }
The second contribution is bounded in a similar way as
    \eqn{
    \sum_{j=1}^{\infty}\expec[X^2 \indic{X\leq m_j}]\eta^{j}
    =\expec\Big[\sum_{j=1}^{\infty} X^2\eta^{j} \indic{X\leq K\eta^{-j}}\Big]
    =\expec\big[\frac{X^2\eta^{b(X)}}{1-\eta}\big],
    }
where $b(x)=\min\{i\colon K\eta^{-j}\geq x\}=\lceil \log{(x/K)}/\log{1/\eta)}\rceil\geq \log{(x/K)}/\log{(1/\eta)}$,
so that $\eta^{b(X)}\leq K/X$. Therefore,
    \eqn{
    \sum_{j=1}^{\infty}\expec[X^2 \indic{X\leq m_j}]\eta^{j}
    \leq \frac{1}{1-\eta} \expec[KX]=\frac{K\nu}{1-\eta},
    }
as required.
\hfill\qed
\medskip

Combining \eqref{E[Z2]} with Lemma \ref{lem-trunc-eff-var} yields:
$$
\lim_{n\to \infty}
	\expec[|\BPbar^{\sss (\vec{m})}(t_n)|^2]
\leq
	\frac{2A^2K\nu}{(1-\eta)},
$$
so that that (b) follows with $C=\frac{2A^2\nu}{(1-\eta)}$ .

For Proposition \ref{prop-sec-mom-clt}(c), we start by showing that, for the original branching process $(|\BP(t)|)_{t\ge 0}$,
	\eqn{
	\label{conv-finvar}
	\e^{-\alpha t}
	\sum_{j=0}^{k_t(x)}
	\expec[|\BP_j[t,t+s)|]
	\to
	A\Phi(x)\FR(s).
	}	
Conditioning on the lifetime (with c.d.f.\
equal to $G$) of the first individual, after which the individual dies and
splits in a random number offspring with mean $\nu$,
	\eqn{
	\label{rec-finvar}
	\expec[|\BP_j[t,t+s)|]=\nu
	\int_0^t \expec[|\BP_{j-1}[t-y,t+s-y)|]\,dG(y).
	}
As before,
	\eqn{
	|\BPbar_j[t,t+s)|=\e^{-\alpha t}|\BP_j[t,t+s)|.
	}
Rewriting (\ref{rec-finvar}) we obtain the recursion
	\eqn{
	\label{ren-rec}
	\expec[|\BPbar_j[t,t+s)|]=
	\int_0^t \expec[|\BPbar_{j-1}[t-y,t+s-y)|]\,d{\bar G}(y).
	}
Hence,  if we continue to iterate, we get
	\eqn{
	\label{conv-rec}
	\expec[|\BPbar_j[t,t+s)|]=
	\int_0^t
	\expec[|\BPbar[t-y,t+s-y)|]\,
	d{\bar G}^{\sss \star j}(y),
	}
where ${\bar G}^{\sss \star j}$ is the $j$-fold convolution of ${\bar G}$, and
hence the distribution function of the independent sum of $j$ copies of a
random variable each having c.d.f. ${\bar G}$. This is the point where \ch{we will use the
CLT.}
For fixed $s>0$, we define
	\eqn{
	\label{def-etan}
	|\BPbar_{\sss >m}[t,t+s)|=
	\sum_{j=m+1}^\infty |\BPbar_j[t,t+s)|.
	}
Observe that $|\BP[t,t+s)|=\sum_{j=1}^\infty |\BP_j[t,t+s)|$ is the total
number of alive individuals of residual lifetime at most
$s$, so that by \cite[Theorem 24.1]{harris1963}, (since $G$ admits a density and $1<\nu<\infty$, the conditions of this theorem are fulfilled),
	$$
	\lim_{t\to\infty}
	\expec[|\BPbar[t,t+s)|]=\lim_{t\to\infty}\sum_{j=0}^\infty\expec[|
	\BPbar_j[t,t+s)|]=A\FR(s),
	$$
where
	\eqn{
	\label{definitie-A}
	A=\frac{\nu-1}{\alpha\nu^2 \int_0^\infty y \e^{-\alpha y}\,dG(y)}=\frac{\nu-1}{\alpha\nu \bar{\nu}}.
	}
Hence, \eqref{conv-finvar} follows if we show that
	\eqn{
	\label{conv-finvar2}
	\expec[|\BPbar_{\sss >k_t(x)}[t,t+s)|]
	\to
	A\FR(s)-A\FR(s)\Phi(x)=A\FR(s)\Phi(-x).
	}
Note that
	\eqn{
	\label{conv-finvar3}
	\expec[|\BPbar_{\sss >k_t(x)}[t,t+s)|]
	=\int_0^t
	\expec[|\BPbar[t-u,t-u+s)|]\,
	d{\bar G}^{\sss \star k_t(y)}(u).
	}
Take an arbitrary $\varepsilon>0$ and take $t_0$ so large so that for $t>t_0$,
	\eqn{
	|\expec[|\BPbar[t,t+s)|]-A\FR(s)|\leq \varepsilon.
	}
Then,
	\eqan{
	\label{conv-finvar4}
	&\Big|\expec[|\BPbar_{\sss >k_t(x)}[t,t+s)|]-A\FR(s)\Phi(-x)\Big|\\
	&\qquad\leq \vep {\bar G}^{\sss \star k_t(x)}(t)
	+A\FR(s)\big|{\bar G}^{\sss \star k_t(x)}(t)-\Phi(-x)\big|+\int_{t-t_0}^{t}
	\expec[|\BPbar[t-u,t-u+s)|]\,
	d{\bar G}^{\sss \star k_t(y)}(u).\nn
	}
The last term vanishes since $\expec[|\BPbar[t,t+s)|]$ is uniformly bounded and ${\bar G}^{\sss \star k_t(y)}(t)-{\bar G}^{\sss \star k_t(y)}(t-t_0)=o(1)$.
Furthermore, with $m=k_t(x)\ra \infty$,
	\eqn{
	\label{identity-t-k}
	k_t(x)\sim\frac{t}{{\bar \nu}}+x\sqrt{t\frac{{\bar \sigma}^2}{{\bar \nu}^3}}
	\qquad\Longleftrightarrow \qquad t\sim m{\bar \nu}-x{\bar \sigma}\sqrt{m}.
	}
As a result, by the CLT and the fact that ${\bar \nu}$ and $\bar{\sigma}^2$
are the mean and the variance of the distribution function $\bar{G}$,
	\eqn{
	\lim_{t\rightarrow \infty} {\bar G}^{\sss \star k_t(x)}(t)=\Phi(-x).
	}
Together with \eqref{conv-finvar4}, this proves the claim in \eqref{conv-finvar2},
\ch{and hence Proposition \ref{prop-sec-mom-clt}(c).}

We continue with the proof of Proposition \ref{prop-sec-mom-clt}(a) for the $n$-dependent 
CTBP. We denote the number of alive  individals at time $t$ in the $n$-dependent CTBP  by $|\BP_{\sss(n)}(t)|$.  We then have to show that \ch{uniformly in $n$,}
	\eqn{
	\label{n-dependent-conv}
	\e^{-\alphan t_n} \expec[|\BP_{\sss(n)}(t_n)|]\to A,
	}
where $A$ is given in (\ref{definitie-A}).
Denote by $\varphi(s)=\int_0^\infty \e^{-sy}g(y)\,dy$, the Laplace transform of the lifetime distribution ($g$ is the density of $G$).
Then
	\eqn{
	\label{Laplace-transform}
	\int_0^\infty \e^{-st} \expec[|\BP_{\sss(n)}(t)|]\,dt =\frac{1-\varphi(s)}{s(1-\nu_n \varphi(s))}.
	}
	\ch{
This equation follows directly from \cite[Equation 16.1]{harris1963}, with $m$ replaced by $\nu_n$ and is valid
when the real part of $s$ satisfies $\Re(s) > \alphan$, where $\alphan>0$ is defined as the unique value  with
$\nu _n \varphi(\alphan)=1$.}  From the inversion formula for
Laplace transforms, we obtain:
	\eqn{
	\expec[|\BP_{\sss(n)}(t)|]=\frac1{2\pi i}\int_{\Gamma} \e^{st}
	\frac{1-\varphi(s)}{s(1-\nu_n \varphi(s))}\,dt,
	}
where $\Gamma$ is the path $(c-i\infty,c+i\infty)$, with $c>\alpha_n$.
Since $\alpha_n\to\alpha$ and $\nu_n\to \nu>1$ and $\varphi(s)$ is the Laplace transform of a {\it probability density}, the function $s(1-\nu_n \varphi(s))$ has a simple zero $s=\alphan$, but no other zeros in a small strip $|s-\alphan|<\varepsilon$. It is now easy to conclude from Cauchy's theorem, calculating the residue at $s=\alpha_n$, \ch{that}
	\begin{eqnarray}
	\expec[|\BP_{\sss(n)}(t_n)|]
	&=&
	\e^{\alphan t_n}\frac{1-\varphi(\alphan)}
	{\alpha_n\cdot(-\nu_n \varphi'(\alpha_n))}\Big(1+O(\e^{-\varepsilon t_n})\Big)\nonumber\\
	&=&
	A_n
	\e^{\alphan t_n}\Big(1+O(n^{-\varepsilon/(2\alphan)})\Big),
	\end{eqnarray}
where
	\eqn{
	\label{definitie-An}
	A_n=\frac{\nu_n-1}{\alphan\nu_n^2 \int_0^\infty y\e^{-\alphan y}\,dG(y)}=\frac{\nu_n-1}{\alphan\nu_n\bar{\nu}_n}.
	}
Since $A_n\to A$, by Condition \ref{cond-degrees-regcond}(b), the claim  \eqref{n-dependent-conv} follows.

For the second statement in Proposition \ref{prop-sec-mom-clt}(a) for the $n$-dependent CTBP, 
we replace \eqref{prod-sum-ineq} 
by the equivalent $n$-dependent statement:
  	\eqn{
	\label{n-depen-ineq}
    	1-\prod_{i=1}^k \frac{\nu_n^{\sss(i)}}{\nu_n}
		\leq (\log(1/\eta))^{-1}
	\expec\big[X_n[\log{(X_n/\Kn)}_+]\big],
    	}
Since $X_n\convd X$, $\expec[X_n]\to \expec[X]$ and $\limsup_{n\rightarrow \infty} \expec\big[X_n \log(X_n/\Kn)_+]\big]=0$ the statement follows.	 

For the $n$-dependent case of Proposition \ref{prop-sec-mom-clt}(b), we need to show that uniformly in $n$,
    	\eqn{
	\label{n-dep-upperb}
    	\e^{-2\alpha_n t_n}\expec\big[|\BP_{\sss (n)}^{\sss (\vec{m})}(t_n)|^2\big]
	\leq C K_n,
	}
for some constant $C$ and where $K$ is defined through $m_i=K_n\eta_n^i$ 
where $\eta_n=\nu_n\int_0^\infty \e^{-2\alpha_n y}\,d G_n(y)$ and $\nu_n=\expec[X_n]$.
Copying the derivation which leads to \ch{ \eqref{driezeventien}}, we obtain:
\eqan{
\label{n-dep-sec-mom}
	\expec[|\BPbar_{\sss(n)}^{\sss (\vec{m})}|^2]
	&=\sum_{j=1}^\infty b_1^{\sss (n)}\cdots b^{\sss (n)}_{j-1}a^{\sss (n)}_j \expec[|{\BPbar}_{\sss(n)}^{\sss (\vec{m}_j)}|]^2*{\bar \GF}^{\sss (n)}_{j-1}*Q^{\sss (n)},
	}
where
	\eqn{
	a^{\sss (n)}_j=\frac{\eta_n \expec[(X_n\wedge m_j)^2 ]}{\nu_n},\qquad
	b^{\sss (n)}_j=\frac{\eta_n  \expec[(X_n\wedge m_j)]}{\nu_n},
	}
	\ch{
	\eqn{
	{\bar \GF}_n(t)= \nu_n\int_0^t \e^{-\alpha_n y}\,d\GF(y),
	\qquad
	Q_n(t)=
	\eta_n^{-1}\nu_n \int_0^t \e^{-2\alpha_n y}\,d\GF(y).
	}
	}
From the proof of Lemma \ref{lem-trunc-eff-var}, we readily obtain that:
	\eqn{
	\label{n-dep-ub}
  	\sum_{j=1}^\infty b_1^{\sss (n)}\ldots b^{\sss (n)}_{j-1}a^{\sss (n)}_j\leq
  	\sum_{j=1}^\infty m_j^2 \prob(X_n>m_j)\eta_n^j+
  	\sum_{j=1}^\infty  \expec[X_n^2 \indic{X_n\leq m_j}]\eta_n^j
	\leq
	2\frac{K_n\nu_n}{1-\eta_n}.
  }
Since, $\nu_n\to\nu$ and \ch{$\eta_n\to \eta$} as $n\ra \infty$, we find, by combining
\eqref{n-dep-sec-mom} with \eqref{n-dep-ub}, that given $\varepsilon >0$, there is an $n_0$ so that for $n>n_0$,

	\eqn{
	\e^{-2\alpha_n t_n}\expec\big[|\BP_{\sss (n)}^{\sss (\vec{m})}(t_n)|^2\big]
	\leq
	\frac{2\Kn(\nu+\varepsilon)}{(1-\eta-\varepsilon){(\bar \nu-\varepsilon)}^2}\leq C\Kn.
	}
By (if necessary) enlarging  the constant $C$ we see that \eqref{n-dep-upperb} holds for all $n$.
This proves Proposition \ref{prop-sec-mom-clt}(b) for the $n$-dependent CTBP.

Finally, we consider Proposition \ref{prop-sec-mom-clt}(c) for the $n$-dependent CTBP. We denote by
$|\BP_{{\sss (n)},j}[t,t+s)|$ the number of individuals in generation $j$ having residual lifetime at most $s$ 
at time $t$ of the CTBP with offspring 
given by $X_n$.  Then,
we obtain, compare \eqref{conv-rec},
	\eqn{
	\expec[|{\BPbar}_{\sss (n),>k}[t,t+s)|]=\int_0^t \expec[|{\BPbar}_{\sss (n)}[t-y,t+s-y)|]
	\, d{\bar G}^{\sss \star k}_{n}(y).
	}
As in the proof of Proposition \ref{prop-sec-mom-clt}(c),
	$$
	\lim_{t\to\infty}
	\expec[|{\BPbar}_{\sss (n)}(t)|]=A_n\FR(s),
	$$
where $A_n$ was defined in \eqref{definitie-An}.
A small extension of the CLT yields that
	$$
	{\bar G}^{\sss (n)}_{k_{t_n}(x)}(t_n)\to
	\Phi(x).
	$$
Since $\nu_n\to \nu$, $\alpha_n\to \alpha$, $A_n\to A$, as $n\to \infty$, it follows that
	\eqn{
	\label{n-dep-lim-untunc-case}
	\expec[|\BPbar_{{\sss (n),>k_{t_n}(x)}}[t,t+s)|]
	=
	\int_0^{t_n} \expec[|\BPbar_{{\sss (n)}}[t_n-y,t_n-y+s)|] \, d{\bar G}^{\sss \star k_{t_n}(x)}_n(y)\to
	A\FR(s)\Phi(x).
	}
This completes the proof of Proposition \ref{prop-sec-mom-clt}(c) for $n$-dependent CTBPs.
\hfill\qed

\section{Coupling to CTBP: Proof of Proposition \ref{prop-CTBP-asy}}
\label{sec-thinning-CM}

\subsection{The coupling}
\label{sec-clus-explor}
The exploration of the total progeny of a branching process satisfies the same recurrence relation as in
\eqref{resurs-S}, apart from the fact
that the random variables $(X_k)_{k\geq 1}$ are i.i.d.\ for a CTBP. For
$\CMnd$, clearly,
$(X_k)_{k\geq 1}$ are not i.i.d. We now describe stochastic relations
between $(X_k)_{k\geq 1}$ given in the CM
and an i.i.d.\ sequence $(Y_k)_{k\geq 1}$ with distribution equal to that
of $D_n^{\star}-1$ for all $k\geq 1$, where
the distribution of $D_n^{\star}-1$ has probability mass function
$(g_k^{\sss (n)})_{k\geq 0}$ defined by
       \eqn{
       \label{def_Dnstar}
       \prob(D_n^{\star}-1=k)=g_k^{\sss (n)}=\frac{k+1}{\ell_n} \sum_{i=1}^n
\indic{d_i=k+1}, \qquad k\ge 0.
       }
Recall that $\Ver_1,\Ver_2$ are chosen uniformly at random from $[n]$.
We continue with the definition of the \emph{size-biased reordering} of $[n]\setminus\{\Ver_1,\Ver_2\}$:

\begin{Definition}[Size-biased reordering]
\label{size-bias-re}
Given the set $[n]=\{1,2,\ldots,n\}$, vertices $\Ver_1, \Ver_2$, and degree sequence $\bfd$, so that
element $i\in[n]$ has degree
$d_i$, a \emph{size-biased reordering} of $[n]\setminus\{\Ver_1,\Ver_2\}$ of size $m_n$ is a random sequence
of (different) elements
$V_1,\ldots,V_{m_n}$, where we select $V_i,\, 1\le i\le m_n$ at
random from the set
$[n]\setminus \{\Ver_1,\Ver_2, V_1,\ldots,V_{i-1}\}$ with probability proportional to
the remaining degrees:
       $$
       \{d_1,\ldots,d_n\}\setminus
       \{d_{\Ver_1},d_{\Ver_2},d_{V_1},\ldots,d_{V_{i-1}}\}.
       $$
\end{Definition}

Let $B_i+1=d_{V_i}$. We let $X_i$ be the number of sibling half-edges 
of $V_i$ that do not create cycles, i.e., are connected to vertices unequal to 
$\{\Ver_1,\Ver_2, V_1, \ldots, V_{i-1}\}$. Thus, clearly, $X_i\leq B_i$.
The above set-up allows us to define the coupling between
$(B_i)_{i\geq 1}$ and $(Y_i)_{i\geq 1}$, where $(Y_i)_{i\geq 1}$, is an
i.i.d. sequence with distribution \eqref{def_Dnstar}.

\begin{Construction}[Coupling of size-biased reordering]
\label{con-coupling-SB-reordering}
We couple $(B_i)_{i\geq 1}$ and $(Y_i)_{i\geq 1}$ in the following way:\\
(a) Draw $Y_i$ as an independent copy of the distribution \ch{in \eqref{def_Dnstar}}. This
can be achieved by
drawing a uniform half-edge $y$ from the total of $\ell_n$ half-edges. Let
$V_i'=V_y$ denote the vertex to which the chosen half-edge is incident, and
let $Y_i=d_{V_i'}-1$.\\
(b) If $V_i'\not\in \{\Ver_1,\Ver_2, V_1, \ldots, V_{i-1}\}$, then
$B_i=Y_i$ and $V_i=V_i'$,
and we say that $V_i$ is \emph{successfully coupled} \ch{with $V_i'$}.\\
(c) If $V_i'\in \{\Ver_1,\Ver_2, V_1, \ldots, V_{i-1}\}$, so that we draw a
half-edge
incident to the set $\{\Ver_1,\Ver_2, V_1, \ldots, V_{i-1}\}$, 
then we \emph{redraw}
a half-edge from the set of half-edges incident to 
$[n]\setminus \{\Ver_1,\Ver_2, V_1, \ldots, V_{i-1}\}$ with
probability proportional to their degree, we let 
$V_i$ denote the vertex incident to the
half-edge drawn, \ch{$B_i=d_{V_i}-1$,} and we say that both $V_i$ and $V_i'$ are
\emph{miscoupled}.\\
(d) We define $X_i$ as the number of the $d_{V_i}-1$ half-edges incident
to vertex $V_i$ that are not paired to a half-edge incident to 
$\{\Ver_1,\Ver_2, V_1, \ldots, V_{i-1}\}$.
\end{Construction}

We next investigate the above coupling. For this, it will be useful to
note that when $D_n$, with distribution function $F_n$ given in
\eqref{def-Fn-CM}, satisfies
Condition \ref{cond-degrees-regcond}(c), then
the maximal degree $\Delta_n=\max_{i\in[n]} d_i$ satisfies
       \eqn{
       \label{Deltan-bd}
       \Delta_n=o(\sqrt{n/\log{n}}).
       }
Indeed, suppose that $\Delta_n\geq \vep \sqrt{n/\log{n}}$. Then, pick
$K_n=n^{1/4}$ to obtain that
       \begin{eqnarray}
       \expec[D_n^2\log{(D_n/\Kn)}_+]&=&\frac1{n}\sum_{k=1}^n d_k^2
\log(d_k/n^{1/4})_+\geq \frac{\Delta_n^2}{n}\log(\Delta_n/n^{1/4})      
       \nn\\
       &\geq & n^{-1} (\vep \sqrt{n/\log{n}})^2 \log(n^{1/4}/(\log{n})^{1/2})
       \geq \vep^2/8.
       \end{eqnarray}
This is in contradiction to Condition \ref{cond-degrees-regcond}(c), so
we conclude that \eqref{Deltan-bd}
holds.

We define the sigma-algebra  ${\mathcal G}_{i}$ by ${\mathcal
G}_{i}=\sigma(d_{\Ver_1}, d_{\Ver_2}, X_0^{\sss(1)}, X_0^{\sss(2)}, X_j,
B_j)_{j\in[i]}$, \ch{see Section \ref{sec-pfs} for the definition of $X_0^{\sss(i)},\, i=1,2$.}

\begin{Lemma}[Coupling to an i.i.d.\ sequence]
\label{lem-coupling-c}
Assume that Condition \ref{cond-degrees-regcond}(c) holds. For all $i\leq
m_n$, and assuming that $m_n\leq \sqrt{n\log{n}}$,
    \eqn{
    \label{coupling_upb}
    \prob(B_i\neq Y_i\mid {\mathcal G}_{i-1})\leq
       \frac{1}{\ell_n} \Big(S_0+\sum_{s=1}^{i-1} (B_s+1)\Big)=o(1),
    }
where, as before, $S_0= X_0^{\sss(1)}+ X_0^{\sss(2)}$.
\end{Lemma}

\proof In  Construction \ref{con-coupling-SB-reordering},
$B_i\neq Y_i$ precisely when $V_i'\in \{\Ver_1,\Ver_2, V_1, \ldots,
V_{i-1}\}$, which, given $ {\mathcal G}_{i-1}$,  has probability at most
    \eqn{
    \prob(B_i\neq Y_i\mid {\mathcal G}_{i-1})
    \leq \frac{1}
    {\ell_n}\Big(S_0+\sum_{s=1}^{i-1} (B_s+1)\Big),
    }
since $Y_i$ draws uniformly from a total of $\ell_n$ half-edges, whereas
in the previous draws at most $X_0^{\sss(1)}+ X_0^{\sss(2)}+\sum_{s=1}^{i-1}
(B_s+1)$ half-edges are incident to the vertices $\{\Ver_1,\Ver_2, V_1, \ldots, V_{i-1}\}$.
  
By \eqref{Deltan-bd}, $\Delta_n=o(\sqrt{n/\log{n}})$ so that
	$$
	X_0^{\sss(1)}+ X_0^{\sss(2)}+\sum_{s=1}^{i-1} (B_s+1)\leq
	m_n\Delta_n=o(n).
	$$
The final statement in \eqref{coupling_upb} follows from $\ell_n= n \mu (1+o(1))$.
\hfill\qed
\medskip

\begin{Lemma}[Probability of drawing a half-edge incident to a previously found vertex]
\label{lem-coupling-b}
Assume that Condition \ref{cond-degrees-regcond}(c) holds.  For all
$i\leq m_n$, and assuming that $m_n\leq \sqrt{n\log{n}}$,
    	\eqn{
	\label{vierzeven}
    	\prob(X_i<B_i\mid {\mathcal G}_{i-1})\leq
       \frac{B_i}{\ell_n(1-o(1))} \Big(S_0+\sum_{s=1}^{i-1} B_s\Big).
    	}
\end{Lemma}

\proof Recall the definition of $S_i$ in \eqref{resurs-S}.
We have $X_i<B_i$ precisely when we pair at least one of the $B_i$
half-edges
to a half-edge incident to $\{\Ver_1,\Ver_2, V_1, \ldots, V_{i-1}\}$.
Since there are precisely $B_i$ half-edges that need to be paired, and
the number of half-edges incident to $\{\Ver_1,\Ver_2, V_1, \ldots, V_{i-1}\}$, 
given ${\mathcal G}_{i-1}$, equals $S_{i-1}$, we find
    	\eqn{
    	\prob(X_i<B_i\mid {\mathcal G}_{i-1})
    	\leq \frac{B_i\cdot S_{i-1}}
    	{\ell_n-\sum_{s=1}^{i-1} (B_s-1)-S_0-1} .
    	}
Clearly, $S_{i-1}\leq S_0+\sum_{s=1}^{i-1} B_s$, which completes the
proof.
As before, $\sum_{s=1}^{i-1} (B_s-1)\leq m_n\Delta_n
=o(n)$ a.s.\ when $m_n\leq \sqrt{n\log{n}}$, \ch{which explains the $\ell_n(1-o(1))$ in the
denominator of \eqref{vierzeven}.}
\hfill\qed
\medskip

\noindent
\paragraph{Coupling the flows.}
In the above, we have described the coupling of our size-biased reordering.
\ch{We now extend this  to a coupling between the SWGs and CTBPs}. We will couple
in such a way that vertices that are succesfully coupled (and thus have
the same number of offspring half-edges, respectively, individuals), 
also receive the \emph{same} weight
along these half-edges. Therefore, the subtrees of successfully coupled 
half-edges in the SWG and individuals in the CTBP are found at precisely
the same times. We will consistently refer to alive objects in the SWG
as alive half-edges, and as alive individuals in the CTBP.

More specifically, offspring half-edges or individuals of miscoupled vertices
are \ch{by definition} miscoupled. \ch{The weights assigned to} miscoupled half-edges in the SWG
and individuals in the CTBP will be independent.
Recall that $T_k$ denotes the time at which the $k$th half-edge is found by
the flow in the SWG. When the half-edge found is incident to
a sucessfully coupled vertex, then we use Construction \ref{con-coupling-SB-reordering}
to couple the number of offspring half-edges
in the SWG to the offspring individuals in the CTBP.
When the half-edge is incident to a miscoupled vertex, it is only present in
the SWG, and we draw a half-edge from the set of available half-edges
as in Construction \ref{con-coupling-SB-reordering}(c),
ignoring $Y_k$ in Construction \ref{con-coupling-SB-reordering}.
\ch{ We define  $(T_k')_{k\ge 0}$ as the times  where an individual dies in the CTBP, but
no half-edge is found by the flow in the SWG. These events result from one or more miscouplings between the CTBP and the SWG. At such times we draw a 
half-edge $y$ uniformly at random from the total number of half-edges,
and let $Y_i=d_{V_y}-1$ denote the number of sibling half-edges. 
Note that in this case, we do not rely on Construction \ref{con-coupling-SB-reordering},
and thus $V_y$ is not part of the size-biased reordering.}

Because of the above construction, differences arising in the coupling are due to two
effects:\\
(1) a \emph{miscoupling} occurs:
miscoupling between the size-biased reordering $B_i$ and the i.i.d.\ draw
from the degree distribution $Y_i$; and\\
(2) a \emph{cycle-creating event} occurs: Here we  refer to the occurrence of
cycles, which makes $X_i<B_i$, and, by our construction of the
collision edges, removes the $B_i-X_i$ half-edges incident to
vertex $V_i$, as well as the $B_i-X_i$ half-edges to which
they are paired from $\SWG$.\\
Recall that offspring of miscoupled vertices are also miscoupled, so any miscoupling 
gives rise to a tress of miscoupled children half-edges in the SWG, respectively, individuals
in the CTBP.

\ch{In order to be able to show convergence in probability to a random variable, we assume that all couplings are defined on one and the same probability space}.

\subsection{Coupling the SWG to a CTBP: proof of Proposition \ref{prop-CTBP-asy}(a)}
\label{sec-perfect-coupling-to-sn}
Consider an age-dependent branching process with lifetimes having a distribution admitting density $g$, and
offspring distribution given by $(f_k)_{k\ge 1}$ in the first generation and offspring
distribution $(g_k)_{k\ge 0}$ in the second and all further generations, \ch{and let $D$ have 
probability mass function (p.m.f.) $(f_k)_{k\ge 1}$.}

%
Let $D^{\star}_n$ ($D^{\star}$) be a random variable such that $D_n^{\star}-1$ ($D^{\star}-1$) has \ch{p.m.f.}
$g_k^{\sss (n)}$ ($g_k$), i.e.,
	\eqn{
	\prob(D_n^{\star}-1=k)=\frac{k+1}{\ell_n} \sum_{i=1}^n \indic{d_i=k+1},\qquad\qquad
	\prob(D^{\star}-1=k)=\frac{(k+1)\prob(D-1=k)}{\mu}.
	}
By Condition \ref{cond-degrees-regcond}, $D^{\star}_n\convd D^{\star}$, so that, since $D^{\star}_n$ and $D^{\star}$ are discrete distributions,
	\begin{equation}
	\label{conv-tv}
	\TVD(D_n^{\star},D^{\star})\to 0,
	\end{equation}
where $\TVD$ denotes the total variation distance, see for instance \cite[Theorem 6.1]{thorisson}.

{\it Proof of Proposition \ref{prop-CTBP-asy}(a).}
Take $s_n$ maximal such that
	\begin{equation}
	\label{speed-coup}
	\e^{2\alpha s_n} \TVD(D_n^{\star},D^{\star})\to 0.
	\end{equation}
According to \eqref{conv-as-BP}, and with $i\in \{1,2\}$,
	\begin{equation}
	\label{as-conv}
	\e^{-\alpha s_n}|\BP^{\sss(i)}(s_n)|\convas\WW^{\sss(i)},
	\end{equation}
where $\WW^{\sss(i)}$ are two independent copies of $\WW$.
Since $\prob(\WW^{\sss(i)}<\infty)=1$ and $\e^{\alpha s_n}\to \infty$, we conclude
that $|\BP(s_n)|\le k_n$ whp if we take $k_n=\lfloor\e^{2\alpha s_n}\rfloor$.
If this $k_n$ does not satisfy $k_n=o(\sqrt{n})$ then lower $s_n$
so that the corresponding value of $k_n=\lfloor\e^{2\alpha s_n}\rfloor$ does satisfy
$k_n=o(\sqrt{n})$.

Recall that ${\mathcal
G}_{i}=\sigma(d_{\Ver_1}, d_{\Ver_2}, X_0^{\sss(1)}, X_0^{\sss(2)}, X_j,
B_j)_{j\in[i]}$. \ch{By Boole's inequality,}
	\eqn{
	\prob(X_i\neq Y_i\mid {\mathcal G}_{i-1})
	\leq \prob(X_i<B_i\mid {\mathcal G}_{i-1})+\prob(B_i\neq Y_i\mid {\mathcal G}_{i-1}).
	}
Consequently, \ch{by Lemmas \ref{lem-coupling-c}--\ref{lem-coupling-b}}, a lower bound for the probability of coupling successfully
during the first \ch{ $k_n=o(\sqrt{n})$} pairings is \ch{
	\begin{eqnarray}
	\label{low-bnd-no-cycle}
	&&\prob\big(\bigcap_{i=1}^{k_n}\{X_i=Y_i\}\big)
	=1-\prob\big(\bigcup_{i=1}^{k_n}\{X_i\neq 	Y_i\}\big)\nonumber\\
	&&\qquad\geq 1-\frac1{n(1-o(1))} \sum_{i=1}^{k_n}\Big(\expec[B_i S_0]+1+\sum_{s=1}^{i-1} (\expec[B_i B_s]+1)\Big)
	\geq 1-c k_n^2/n\to 1.
	\end{eqnarray}
	}
Here we rely on the fact that \ch{ 
	\eqn{
	\expec[B_i\mid {\mathcal G}_{i-1}] 
	\leq \sum_{j\in [n]} \frac{d_j(d_j-1)}{\ell_n-2k_n\Delta_n}=\nu_n(1+o(1)),
	}
whenever $k_n\Delta_n=o(n)$, which follows from \eqref{Deltan-bd}.}

The lower bound \eqref{low-bnd-no-cycle} implies that, whp,
the shortest weight graph $(\SWG(s))_{s\le s_n}$ is perfectly 
coupled to the CTBP $(\BP(s))_{s\leq s_n}$.
This proves Proposition \ref{prop-CTBP-asy}(a). \hfill\qed

We close this section by investigating moments of the size-biased variables
$(B_i)_{i\geq 1}$ arising in the size-biased reordering.
These moments play a crucial role throughout the remainder of this paper, and allow
us to compare $(B_i)_{i\geq 1}$ to an i.i.d.\ sample of random variables having the 
size-biased random distribution.

\begin{Lemma}[Moments of the size-biased reordering]
\label{lem-mom-size-biased}
Assume that Condition \ref{cond-degrees-regcond}(a-c) holds.  For all
$i\leq m_n$, and assuming that $m_n\leq\sqrt{n\log{n}}$, and for any 
$\Kn\ra \infty$ such that $\Kn^2=o(n/m_n)$, 
    	\eqan{
    	\expec[B_i\indic{B_i\leq \Kn}\mid {\mathcal G}_{i-1}]&=(1+\op(1))\nu_n,
	\label{size-biased-first}\\
    	\expec[B_i\indic{B_i>\Kn}\mid {\mathcal G}_{i-1}]&=\op(1).
	\label{size-biased-trunc-first}
    	}
\end{Lemma}

\proof \ch{We use the upper bound}
	\eqn{
	\expec[B_i\indic{B_i\geq a}\mid {\mathcal G}_{i-1}]\leq \frac{1}{\ell_n(1-o(1))} 
	\sum_{l\in [n]} d_l(d_l-1)\indic{d_l\geq a+1},
	}
where we again use that, since $m_n\leq \sqrt{n\log{n}}$, 
	\eqn{
	\ell_n-S_0-\sum_{j=1}^{i-1} B_j\geq \ell_n-m_n\Delta_n=\ell_n-o(n).
	}
This provides the necessary upper bound in \eqref{size-biased-first} by taking $a=0$ and
\ch{from the identity} $\nu_n=\sum_{l\in [n]} d_l(d_l-1)/\ell_n$. For \eqref{size-biased-trunc-first},
this also proves the necessary bound, since 
	\eqn{
	\frac1{\ell_n}
	\sum_{l\in [n]} d_l(d_l-1)\indic{d_l\geq \Kn}=o(1).
	}
For the lower bound in \eqref{size-biased-first}, we bound, instead,
	\eqn{
	\expec[B_i\indic{B_i\leq \Kn}\mid {\mathcal G}_{i-1}]
	\geq \frac{1}{\ell_n} \sum_{l\in [n]} d_l(d_l-1)\indic{d_l\leq \Kn+1}\indic{l \text{ not chosen yet}},
	}
where the event  ``$l$ is not chosen yet"  means that the vertex $l$ has not been chosen in
the size-biased reordering until time $i-1$. We now bound
	\eqn{
	\expec[B_i\indic{B_i\leq \Kn}\mid {\mathcal G}_{i-1}]
	\geq \frac{1}{\ell_n} \sum_{l\in [n]} d_l(d_l-1)\indic{d_l\leq \Kn+1}
	-\frac{1}{\ell_n} \sum_{l\in [n]} d_l(d_l-1)\indic{d_l\leq \Kn+1}\indic{l \text{ is chosen}}.
	}
The first term equals $\nu_n(1+o(1))$. The second term is a.s.\ bounded by $m_n \Kn^2/\ell_n=o(1)$,
since $\Kn^2=o(n/m_n)$. \hfill\qed

\subsection{Completing the coupling: Proof of Proposition \ref{prop-CTBP-asy}(b)}
\label{sec-coupling-completed}
In this section, we use Proposition \ref{prop-sec-mom-clt} to prove Proposition \ref{prop-CTBP-asy}(b).
In order to bound the difference between $\BP(t)$ and $\SWG(t)$, we will introduce several events.
Let $\Bn, \Cn, \vep_n, \mnup, \mnlow$ denote sequences of constants for which $\Bn, \Cn\ra \infty$ and 
$\vep_n\ra 0$ arbitrarily slowly, and $\mnup\gg \sqrt{n}, \mnlow\ll\sqrt{n}$. Later in this proof,
we will formulate restrictions on these sequences.

Define the event $\CA_n$ as follows:
	\eqan{
	\label{symm-diff-small}
	\CA_n&=\{|\SWG(t_n+\Bn)\symdiff \BP_{\sss(n)}(t_n+\Bn)|< \vep_n\sqrt{n}\},
	}
where \ch{$|\SWG(t_n+\Bn)\symdiff \BP_{\sss(n)}(t_n+\Bn)|$} is the number of miscoupled half-edges
plus the number of miscoupled individuals.
Then Proposition \ref{prop-CTBP-asy}(b) can be reformulated as 
	\eqn{
	\label{aim-An-event}
	\prob(\CA_n^c\mid \FF_{s_n})=\op(1).
	}

In order to prove \eqref{aim-An-event},
we introduce the following events:
	\eqan{
	\BB_n&=\{\YBP(t_n+\Bn)\leq \mnup\} \cap \{\YSWG(t_n+\Bn)\leq \mnup\}\nn\\
	&\qquad \qquad \cap \{\YBP(t_n-\Bn)\leq \mnlow\} \cap \{\YSWG(t_n-\Bn)\leq \mnlow\},
	\label{Y-bd}\\
	\CC_n&=\{\SWG(t)=\BP(t),~\forall t\leq t_n-\Bn\},
	\label{perfect-coupling-below-tn}\\
	\!\!\DD_n&=\{\not\exists i \text{ such that } T_i\leq t_n+\Bn, V_i\text{ miscoupled }, d_{V_i} \geq \Cn\},
	\label{perfect-coupling-tn-small}
	}
where 
	\eqn{
	\YBP(t)=|\{v\colon v\in \BP_{\sss(n)}(s) \text{ for some }s\leq t\}|,
	}
denotes the total number of individuals ever born into the BP before time $t$ and
	\eqn{
	\YSWG(t)=|\{v\colon v\in \SWG(s) \text{ for some }s\leq t\}|,
	}
denotes the number of half-edges in the SWG that have ever been alive before time $t$.
Informally, on $\BB_n$, the total number of half-edges in SWG and individuals in the CTBP
are not too large. On $\CC_n$ there is no early miscoupling, while on $\DD_n$, 
there is no miscoupled vertex having high degree until a late stage.

\ch{Obviously
	\eqan{
	&\prob(\CA_n^c\mid \FF_{s_n})\\
	&\qquad\leq \prob(\BB_n^c\mid \FF_{s_n})+\prob(\CC_n^c\cap \BB_n\mid \FF_{s_n})
	+\prob(\DD_n^c\cap \BB_n\cap \CC_n\mid \FF_{s_n})
	+\prob(\CA_n^c\cap\BB_n\cap \CC_n\cap \DD_n\mid \FF_{s_n}).\nn
	}
	}
To bound conditional probabilitites of the form 
\ch{$\prob({\mathcal E^c}\mid \FF_{s_n})$}, we start by noting that it suffices to
prove that $\prob({\mathcal E}^c)=o(1)$, since then, by the Markov \ch{inequality}
and for every $\vep>0$,
	\eqn{
	\prob\Big(\prob({\mathcal E}^c\mid \FF_{s_n})\geq \vep\Big)
	\leq \expec[\prob({\mathcal E}^c\mid \FF_{s_n})]/\vep =\prob({\mathcal E}^c)/\vep=o(1).
	}
Thus, we are left to prove that
	\eqn{
	\prob(\BB_n^c)=o(1),
	\quad
	\prob(\CC_n^c\cap\BB_n)=o(1),
	\quad
	\prob(\DD_n^c\cap\BB_n\cap \CC_n)=o(1),
	\quad
	\prob(\CA_n^c\cap\BB_n\cap \CC_n\cap \DD_n)=o(1).
	}
We will do so in the above order.
	
\begin{Lemma}[Expected number of particles born]
\label{lem-particles-born}
For all $t\geq 0$, \ch{
	\eqn{
	\label{expecY(t)-eq}
	\expec[\YBP(t)]=
	2\Big(1-\frac{\mu_n G(t)}{\nu_n-1}\Big)+
	\frac{\nu_n}{\nu_n-1}\expec[|\BP_{\sss(n)}(t)|].
	}
	}
Consequently, when $\e^{\alphan(t_n+\Bn)}=o(\mnup), \e^{\alphan(t_n-\Bn)}=o(\mnlow)$,
	\eqn{
	\prob(\BB_n^c)=o(1).
	}
	
\end{Lemma}

\proof \ch{Note that we grow two SWGs and two BPs, which explains the factor $2$ in \eqref{expecY(t)-eq}}. 
\ch{
As is well known the expected number of descendants in generation $k$ of a BP equals  $\nu_n^{k}$, where $\nu_n$ denotes the mean offspring. Here, we deal with a delayed BP where in the first generation the mean number of offspring equals $\mu_n=\expec[D_n]$; the factor
$G^{\sss \star k}(t)-G^{\sss \star (k+1)}(t)$ represents the probability that an individual of generation $k$ is alive at time $t$; together this yields:}
	\eqn{
	\expec[|\BP_{\sss (n)}(t)|]=\sum_{k=1}^{\infty} 2\mu_n\nu_n^{k-1} [G^{\sss \star k}(t)-G^{\sss \star (k+1)}(t)],
	\qquad
	\expec[\YBP(t)]=2+\sum_{k=1}^{\infty} 2\mu_n \nu_n^{k-1} G^{\sss \star k}(t).
	}
We can rewrite the equality for $\expec[|\BP_{\sss(n)}(t)|]$ to obtain
	\eqan{
	\expec[|\BP_{\sss(n)}(t)|]&=2\mu_n G(t)/\nu_n+\sum_{k=1}^{\infty} 2\mu_n [\nu_n^{k-1}-\nu_n^{k-2}]G^{\sss \star k}(t)\\
	&=2\mu_n G(t)/\nu_n-2(1-\nu_n^{-1})+(1-\nu_n^{-1})\expec[\YBP(t)].\nn
	}
Solving for $\expec[\YBP(t)]$ yields the proof of \eqref{expecY(t)-eq}.

To bound $\prob(\BB_n^c)$,
we note that we have to bound events of the form $\prob(\YBP(t)\geq m)$
and $\prob(\YSWG(t)\geq m)$ for various choices of $m$ and $t$.
We use the \ch{Markov inequality} and \eqref{expecY(t)-eq} to bound
	\eqn{
	\prob(\YBP(t)\geq m)\leq \expec[\YBP(t)]/m
	\leq \frac{\nu_n}{m(\nu_n-1)}\expec[|\BP_{\sss (n)}(t)|]+\frac{2}{m}.
	}
By Proposition \ref{prop-sec-mom-clt}(d), \ch{$\expec[|\BP_{\sss (n)}(t_n)|]=A_n\e^{\alphan t_n}(1+o(1))$},
so that
	\eqn{
	\prob(\YBP(t_n)\geq m)= \Theta(\e^{\alphan t_n}/m).
	}
The conditions on $t$ and $m$ in Lemma \ref{lem-particles-born} have been chosen precisely
so that $\e^{\alphan (t_n-B_n)}/\mnlow\ra 0$, and $\e^{\alphan (t_n+B_n)}/\mnup\ra 0$.

We continue with $\prob(\YSWG(t)\geq m)$. We use the same steps as above,
and start by computing
	\eqn{
	\expec[\YSWG(t)]=2+2\sum_{k=0}^{\infty} G^{\sss \star k}(t) \expec[P_k^{\star}],
	}
where $P_0^{\star}=\ell_n/n$ and
$$
P_k^{\star}=\sum_{|\pi|=k, \pi\subseteq \CMnd} (d_{\pi_k}-1)/n,
$$ 
is the sum of the number of  half-edges at the ends of paths of lengths $k$ in $\CMnd$,
from a uniformly selected starting point.
See  \cite[Section 5]{Jans09b} for more details on paths in $\CMnd$.
We compute that
\ch{
	\eqn{
	\expec[P_k^{\star}]=\frac{1}{n}\sum_{v_0, \ldots, v_k} d_{v_0} \prod_{i=1}^k \frac{d_{v_i}(d_{v_i}-1)}{\ell_n-2i+1},
	}
	}
where the sum is over distinct vertices in $[n]$. By \cite[Proof of Lemma 5.1]{Jans09b}, the latter sum is bounded 
by
	\eqn{
	\expec[P_k^{\star}]\leq \expec[D_n] \nu_n^{k}/n.
	}
As a result, we have that $\expec[\YSWG(t)]\leq \expec[\YBP(t)]$, and we can repeat our arguments
for $\expec[\YBP(t)]$.
\hfill\qed

\begin{Lemma}[No early miscoupling]
\label{lem-no-early-miscoupling}
When $\e^{\alphan(t_n-\Bn)}=o(\mnlow)$ and $\mnlow=o(\sqrt{n})$, then:
\ch{
	\eqn{
	\prob(\CC_n^c\cap \BB_n)=o(1).
	}
	}
\end{Lemma}

\proof By \ch{the proof of} Lemma \ref{lem-particles-born}, whp $\YBP(t_n-\Bn)\leq \mnlow$.
By \eqref{low-bnd-no-cycle}, the probability that there exists a miscoupling before the
draw of the $\mnlow$th half-edge is $o(1)$ when $\mnlow=o(\sqrt{n})$.
\hfill\qed

\begin{Lemma}[No late miscouplings of high degree]
\label{lem-no-late-miscouplings-high-degree}
If $\mnup\leq\sqrt{n\log{n}}$, and
$\Cn$ satisfies
	\eqn{
	\frac{\mnup^2}{\ell_n}\sum_{i\in[n]} d_i^2 \indic{d_i\geq \Cn}=o(n),
	}
then
	\eqn{
	\prob(\DD_n^c\cap \BB_n\cap \CC_n)=o(1).
	}
\end{Lemma}

\proof \ch{On the event $\BB_n$}: $\YBP(t_n+\Bn)\leq \mnup$.
An upper bound for the probability of miscoupling a vertex of degree at least $\Cn$
during the first $\mnup$ pairings is thus \ch{
	\eqn{
	\label{up-bnd-cycle}
	\sum_{i\in[n]} \Big(\frac{d_i}{\ell_n(1-o(1))}\Big)^2 \indic{d_i\geq \Cn}
	\leq
	\frac1{\ell_n^2(1-o(1))} \sum_{i\in [n]} d_i^2 \indic{d_i\geq \Cn}=o(1).
	}
	}
\hfill\qed

\begin{Proposition}[Miscoupled vertices have small offspring]
\label{prop-no-late-miscouplings-high-degree}
If $\mnup\leq \sqrt{n\log{n}}$ and\\
$\e^{2\alphan \Bn}\Cn \mnup^2/\ell_n=o(\sqrt{n})$, 
then
	\eqn{
	\prob(\CA_n^c\cap \BB_n\cap \CC_n\cap \DD_n)=o(1).
	}
\end{Proposition}

\proof We split the proof into three contributions, namely, a bound on
$|\SWG(t)\setminus \BP_{\sss (n)}(t)|$, a bound on the contribution to $|\BP_{\sss (n)}(t)\setminus \SWG(t)|$
due to cycle-creating events, and a bound on the contribution to $|\BP_{\sss (n)}(t)\setminus \SWG(t)|$
due to miscouplings. We start with the first bound:

\paragraph{A bound on $|\SWG(t)\setminus \BP_{\sss (n)}(t)|$.} 
By construction, the number of miscoupled half-edges in $\SWG(t)$ at any time $t$ is bounded
from above by
	\eqn{
	\sum_{j=1}^{\Mis(t)} \YSWG_j(t-\tilde{T}_j),
	}
where $\Mis(t)$ denotes the number of miscoupled vertices, $\tilde{T}_j$ is the birth
of the $j$th miscoupled vertex, and for the $j$th miscoupled vertex
$\tilde{V}_j$, $\YSWG_j(t)$ is the number of half-edges at flow distance (total edge weight) at 
most $t$ from $\tilde{V}_j$.
On the event $\CC_n$, $\tilde{T}_1\geq t_n-\Bn$. Therefore, for every $t\leq t_n+\Bn$,
on the event $\CC_n$,
	\eqn{
	|\SWG(t)\setminus \BP_{\sss (n)}(t)|\leq \sum_{j=1}^{\Mis(t_n+\Bn)} \YSWG_j(2\Bn).
	}
By the Markov inequality,
	\eqan{
	&\prob\Big(\{|\SWG(t)\setminus \BP_{\sss (n)}(t)|\geq \vep\sqrt{n}\}\cap \BB_n\cap \CC_n\cap \DD_n\Big)\nn\\
	&\qquad \leq (\vep\sqrt{n})^{-1} \expec\Big[\indicwo{\BB_n\cap \CC_n\cap \DD_n}
	\sum_{j=1}^{\Mis(t_n+\Bn)} \YSWG_j(2\Bn)\Big].
	\label{Markov-bd-miscouplings}
	}
We rewrite 
	\eqan{
	&\expec\Big[\indicwo{\BB_n\cap \CC_n\cap \DD_n}
	\sum_{j=1}^{\Mis(t_n+\Bn)} \YSWG_j(2\Bn)\Big]\\
	&\qquad=(1+o(1))\sum_{i\in [n]} \prob(i \text { miscoupled}, \BB_n\cap \CC_n\cap \DD_n)
	\expec[ \YSWG(2\Bn)]\nn\\
	&\qquad\leq \sum_{i\in [n]}\Big(\frac{d_i\mnup}{\ell_n}\Big)^2 \expec[ \YSWG(2\Bn)],\nn
	}
where we use that, upon miscoupling of vertex $i$, we redraw a vertex from the size-biased distribution,
for which the number of half-edges found before time $2\Bn$ is equal to $\expec[ \YSWG(2\Bn)](1+o(1))$
since $\mnup\leq \sqrt{n\log{n}}$ and $\BB_n$ occurs. 
Since $\expec[\YSWG(t)]\leq \expec[\YBP(t)]$, we obtain that
	\eqn{
	\expec[ \YSWG_j(2\Bn)]\leq A_n\e^{2\alphan \Bn}(1+o(1)).
	}
Therefore, we arrive at
	\eqn{
	\label{bd-expec-miscoupling}
	\expec\Big[\indicwo{\BB_n\cap \CC_n\cap \DD_n}
	\sum_{j=1}^{\Mis(t_n+\Bn)} \YSWG_j(2\Bn)\Big]
	\leq O(1) \e^{2\alphan \Bn}\mnup^2/\ell_n.
	}
 Combining \eqref{Markov-bd-miscouplings}-\eqref{bd-expec-miscoupling}
proves that $|\SWG(t)\setminus \BP_{\sss (n)}(t)|=\op(\sqrt{n})$, since 
$ (\vep\sqrt{n})^{-1}\mnup^2/\ell_n=o(1)$.

\paragraph{Bounding the contribution to $|\BP_{\sss (n)}(t)\setminus \SWG(t)|$ due to cycle-creating events.}
On the event $\CC_n$, $\tilde{T}_1\geq t_n-\Bn$. Recall that the $j$th miscoupled vertex
is denoted by $\tilde{V}_j$, and that the time 
of the occurrence of the $j$th miscoupled vertex is $\tilde{T}_j$.
On the event $\DD_n$,  $d_{\tilde{V}_j}\leq \Cn$ for every $j$ for which $\tilde{T}_j\leq t_n+\Bn$.
When a cycle-creating event occurs, the two half-edges that form the 
last edge in the cycle are removed from $\SWG(t)$, but they are kept in $\BP(t)$.
On the event $\BB_n\cap \CC_n$, the expected number of cycle-creating events is bounded 
by
	\eqn{
	O\Big(\frac{\mnup^2}{\ell_n}\Big).
	}
Furthermore, on the event $\BB_n\cap \CC_n$, the expected offspring of the 
half-edges involved in cycle-creating events is at most
	\eqn{
	O\Big(\frac{\mnup^2}{\ell_n}\Big)\expec[\tilde{Y}(2\Bn)],
	}
where $(\tilde{Y}(t))_{t\geq 0}$ is the total number of individuals
that have ever been alive in a CTBP where all individuals have 
i.i.d.\ offpring with law $D_n^{\star}-1$, starting from $D_n^{\star}-1$
individuals. Indeed, we have no information
about the remaining lifetime of the half-edge involved in the 
cycle-creating event. As a result, we \emph{instantaneously} pair it
to an i.i.d.\ draw of a half-edge, and start the $\BP_{\sss (n)}(t)$ 
from there. The total number of individuals ever alive only increases 
by this change. \ch{On the event $\DD_n$, we have:
$\expec[\tilde{Y}(2\Bn)]\leq C_n A_n\e^{2\alphan \Bn}(1+o(1))$.
By assumption $\mnup^2C_n A_n\e^{2\alphan \Bn}/\ell_n=o(\sqrt{n})$.}
Therefore, the contribution to $|\BP_{\sss (n)}(t)\setminus \SWG(t)|$ due to cycle-creating events is $\op(\sqrt{n})$, as required. 

\paragraph{Bounding the contribution to $|\BP_{\sss (n)}(t)\setminus \SWG(t)|$ due to miscouplings.}
We complete the proof of Proposition \ref{prop-no-late-miscouplings-high-degree}
by dealing with the contribution to $|\BP_{\sss (n)}(t)\setminus \SWG(t)|$ due to miscouplings. 
Now, for $\BP_{\sss (n)}(t)$, we do \emph{not} redraw the random
variable $\tilde{V}_j$. We can give an upper bound on the contribution to $\BP_{\sss (n)}(t)$ 
of these miscouplings by instantaneously pairing the half-edges to an i.i.d.\
draw of a half-edge. As a result, the contribution to $|\BP_{\sss (n)}(t)\setminus \SWG(t)|$ 
due to miscouplings can be bounded above by
	\eqn{
	\sum_{j\colon \tilde{T}_j\leq t_n+\Bn} \YBP_j(t_n+\Bn-\tilde{T}_j).
	}
We use that, on the event $\CC_n$, $\tilde{T}_j\geq t_n-\Bn$,
and on the event $\BB_n$, the expected number of miscoupling is at most
$O(\mnup^2/\ell_n)$. Finally, on the event $\DD_n$, $d_{\tilde{V}_j}\leq \Cn$ 
for every $j$ for which $\tilde{T}_j\leq t_n+\Bn$. Therefore,
	\eqn{
	\expec\Big[\YBP_j(t_n+\Bn-\tilde{T}_j)\indicwo{\BB_n\cap \CC_n\cap \DD_n}\Big]
	\leq \Cn A_n\e^{2\alphan \Bn}(1+o(1)).
	} 
We conclude that, on the event 
$\BB_n\cap \CC_n\cap \DD_n$, 
	\eqn{
	\expec\Big[\sum_{j\colon \tilde{T}_j\leq t_n+\Bn} \YBP_j(t_n+\Bn-\tilde{T}_j)\indicwo{\BB_n\cap \CC_n\cap \DD_n}\Big]
	\leq O\Big( \frac{\mnup^2}{\ell_n}\Big)\Cn\e^{2\alphan\Bn}.
	}
By assumption the r.h.s.\ is $o(\sqrt{n})$.
Therefore, the contribution to $|\BP(t)\setminus \SWG(t)|$ due to miscouplings
is $\op(\sqrt{n})$, as required. This completes the proof of 
Proposition \ref{prop-no-late-miscouplings-high-degree}.
\hfill\qed
\medskip

\noindent
{\it Proof of Proposition \ref{prop-CTBP-asy}(b).} It suffices to show that we can choose
the sequences $\Bn,\Cn,\vep_n,\mnup,\mnlow$ such that all conditions in 
Lemmas \ref{lem-particles-born}--\ref{lem-no-late-miscouplings-high-degree} 
and Proposition \ref{prop-no-late-miscouplings-high-degree}  apply. \ch{It is readily verified that we can take:	
	\eqn{
	\label{mn-choices}
	\mnlow=\sqrt{n}/(\log\log{n})^{\alpha/2}, 
	\qquad 
	\mnup=\sqrt{n}(\log{n})^{1/4},
	}
and
	\eqn{
	\label{Bn-Cn-def}
	\Bn=\log\log\log{n},
	\qquad
	\Cn=n^{1/4},
	\qquad
	\vep_n=1/\log{n}.
	}
	}
By Condition \ref{cond-degrees-regcond}(c)
\ch{
	\eqn{
	\frac{1}{n}\sum_{i\in[n]} d_i^2 \indic{d_i\geq \Cn}=\expec[D_n^2\indic{D_n\geq \Cn}]
	\leq \frac{\expec\Big[D_n^2(\log(\e\cdot D_n/\Cn))_+\Big]}{\log n}=o((\log n)^{-1}),
	}
	}
all the conditions in Lemmas \ref{lem-particles-born}--\ref{lem-no-late-miscouplings-high-degree} 
as well as Proposition \ref{prop-no-late-miscouplings-high-degree} hold. Therefore, the claim 
in \eqref{aim-An-event} follows, which completes the proof of Proposition \ref{prop-CTBP-asy}(b).
\hfill\qed


\section{Height CLT and stable age: Proof of Proposition \ref{prop-CLT-stable-age}}
\label{sec-CTBP-pf-strong}

\ch{We first prove Proposition \ref{prop-CLT-stable-age}(a).} By Proposition \ref{prop-CTBP-asy}(a), 
at time $s_n$, \whp,~$(\SWG(s))_{s\leq s_n}$ is perfectly coupled to the two independent 
CTBPs \ch{$(\BP(s))_{s\leq s_n}$}. \ch{The proof contains several key steps.}

\paragraph{\ch{Reduction to a single BP.}} We start by showing that, 
in order for Proposition \ref{prop-CLT-stable-age}(a) to hold, it suffices to 
prove that for $j\in \{1,2\}$, $x,t\in \Rbold$ and $s>0$,
	\eqan{
	\label{aim-clt}
	&\e^{-\alphan t_n}|\BP^{\sss(j)}_{\sss \leq k_{t_n}(x)}[\tn+t,\tn+t+s)|
	\convp \e^{\alpha t} \Phi(x)\FR(s)\sqrt{\WW^{\sss(j)}/\WW^{\sss(3-j)}},
	}
where we use \eqref{WW-vep} in Proposition \ref{prop-CTBP-asy}(a)
to see that $\sqrt{\WW^{\sss(j)}/\WW^{\sss(3-j)}}\in[\vep,1/\vep]$ whp.
Indeed, by Proposition \ref{prop-CTBP-asy}(b) and the fact that $\e^{-\alphan t_n}=n^{-1/2}$, 
\eqref{aim-clt} immediately implies that
	\eqan{
	\label{aim-clt2}
	\e^{-\alphan t_n}|\SWG^{\sss(j)}_{\sss \leq k_{t_n}(x)}[\tn+t,\tn+t+s)|
	&=\e^{-\alphan t_n}|\BP^{\sss(j)}_{\sss \leq k_{t_n}(x)}[\tn+t,\tn+t+s)|+\e^{-\alphan t_n}\op(\vep_n\sqrt{n})\nn\\
	&\convp \e^{\alpha t} \Phi(x)\FR(s)\sqrt{\WW^{\sss(j)}/\WW^{\sss(3-j)}},
	}
Therefore, independence and \eqref{aim-clt} also proves  \ch{that for $j\in \{1,2\}$, $x,y,t\in \Rbold$ and $s_1,s_2$,}
   	\eqan{
	\label{aim-clt-full-SWG}
    	&\e^{-2\alphan t_n}
	|\SWG_{\sss \leq k_{t_n}(x)}^{\sss(j)}[\tn+t,\tn+t+s_1)|
	|\SWG_{\sss \leq k_{t_n}(y)}^{\sss(3-j)}[\tn+t,\tn+t+s_2)|\\
	&\qquad\convp \e^{2\alpha t} \Phi(x)\Phi(y)\FR(s_1)\FR(s_2),\nn
    	}
which is the statement in Proposition \ref{prop-CLT-stable-age}(a).

\paragraph{Using the branching property.}
To prove \eqref{aim-clt}, we note that $(\BP^{\sss(j)}(s))_{s\geq s_n}$ is
the collection of alive individuals in the different generations of a CTBP, starting from the alive particles in
$(\BP^{\sss(j)}(s_n))$. We condition on
$(\BP(s)^{\sss(j)})_{s\in [0,s_n]}$. Then
	\eqn{
	|\BP^{\sss(j)}_{\sss \leq k_{t_n}(x)}[\tn+t,\tn+t+s)|
	=\sum_{i\in \BP^{\sss(j)}(s_n)} \sum_{k=1}^{k_{t_n}(x)-G^{\sss(j)}_i}
	|\BP_k^{\sss(i,j)}[\tn+t-s_n-R_i,\tn+t+s-s_n-R_i)|,
	}
where $G^{\sss(j)}_i$ is the generation of $i\in \BP^{\sss(j)}(s_n)$,
while $R_i$ is its remaining lifetime, and $(\BP^{\sss(i,j)}(t))_{t\geq 0}$ are
i.i.d.\ CTBPs for different $i$, for which the offspring 
for each individual has distribution $D^{\star}_n$.

\paragraph{\ch{Truncating the branching process.}}
We continue by proving that we can truncate the branching process at the expense of 
an error term that converges to zero in probability. 
We let $\BP^{\sss(i,j,\vec{m})}$ denote the branching process $\BP^{\sss(i,j)}$ 
obtained by truncating particles in
generation $l$ (measured from the root $i$) by $m_l=\Kn \eta^{-l}$.
We take $\Kn\rightarrow \infty$ such that $\Kn\e^{-\alphan s_n}=o(1)$.
We first show that, as $t\rightarrow \infty$, we can replace
$\e^{-\alphan t_n}|\BP_{\sss \leq k_{t_n}(x)}^{\sss(i,j)}[\tn,\tn+s)|$ by $\e^{-\alphan t_n}|\BP_{\sss \leq k_{t_n}(x)}^{\sss(i,j,\vec{m})}[\tn,\tn+s)|$, at the expense of a $\op(1)$-term. Indeed,  with
	\eqn{
	|\BP^{\sss(i,j)}(t)|=\sum_{k=1}^{\infty} |\BP_k^{\sss(i,j)}(t)|,
	\qquad
	|\BP^{\sss(i,j,\vec{m})}(t)|=\sum_{k=1}^{\infty}
	|\BP_k^{\sss(i,j, \vec{m})}(t)|,
	}
we have, uniformly in $t\geq 0$ and $k\geq 0$, by the $n$-dependent version of 
Proposition \ref{prop-sec-mom-clt}(a) in Proposition \ref{prop-sec-mom-clt}(d),
	\eqan{
	&\e^{-\alphan t}\expec\Big[|\BP_{\sss \leq k}^{\sss(i,j)}[t,t+s)|-
	|\BP_{\sss \leq k}^{\sss(i,j,\vec{m})}[t,t+s)|\Big] \leq
	\e^{-\alphan t}\expec\Big[
	|\BP^{\sss(i,j)}(t)|
	-|\BP^{\sss(i,j,\vec{m})}(t)|\Big]=o(1).\nn
	}
Therefore, using that the law of $\BP^{\sss(i,j)}_{\sss \leq k_{t_n}(x)}$ only depends on $\FF_{s_n}$
through $R_i, \tn$,
	\eqan{
	&\e^{-\alphan t_n}\sum_{i\in \BP^{\sss(j)}(s_n)} \sum_{k=1}^{k_{t_n}(x)-G^{\sss(j)}_i}\expec\Big[
	|\BP^{\sss(i,j)}_{k}(\tn+t-s_n-R_i)|
	-|\BP^{\sss(i,j,\vec{m})}_{k}(\tn+t-s_n-R_i)|
	\mid \FF_{s_n}\Big]\nn\\
	&\qquad \leq\e^{-\alphan t_n}\sum_{i\in \BP^{\sss(j)}(s_n)} \expec\Big[
	|\BP^{\sss(i,j)}_{\sss \leq k_{t_n}(x)}(\tn+t-s_n-R_i)|
	-|\BP^{\sss(i,j,\vec{m})}_{\sss \leq k_{t_n}(x)}(\tn+t-s_n-R_i)|
	\mid R_i, \tn\Big]\nn\\
	&\qquad =o(1) \sum_{i\in \BP^{\sss(j)}(s_n)}
	\e^{\alphan (\tn-t_n+t-s_n-R_i)}=\op(1) \e^{-\alpha s_n}\sum_{i\in \BP^{\sss(j)}(s_n)}
	\e^{-\alpha R_i},
	}
since the random variable $|\tn-t_n|$ is tight, and
assuming that $s_n\ra\infty$ so slowly that $s_n|\alphan-\alpha|=o(1)$. 
By \cite{jagers1975branching, jagers-nerman4} and with 
$\sigma_i=s-R_i$, the birth-time of individual $i$, 
	\eqn{
	\label{Ms-convp}
	M_{s_n}^{\sss(j)}=\e^{-\alphan s_n}\sum_{i\in \BP^{\sss(j)}(s_n)}
	\e^{-\alphan R_i}\convp \WW^{\sss(j)}/A,
	}
where we use the fact that $(M_s^{\sss(j)})_{s\ge 0}$ is an $n$-independent martingale
by the remark on \cite[p. 234 line 7]{jagers1984growth}. We conclude that
	\eqan{
	\label{sum-trunc}
	&|\BP^{\sss(j)}_{\sss \leq k_{t_n}(x)}[\tn+t,\tn+t+s)|\\
	&\qquad=\sum_{i\in \BP^{\sss(j)}(s_n)} \sum_{k=1}^{k_{t_n}(x)-G^{\sss(j)}_i}
	|\BP_k^{\sss(i,j, \vec{m})}[\tn+t-s_n-R_i,\tn+t+s-s_n-R_i)| +\op(1).\nn
	}

\paragraph{\ch{A conditional second moment method: first moment.}}
We next use a conditional second moment estimate on the sum on the right-hand side 
of \eqref{sum-trunc}, conditionally on $\FF_{s_n}$.
By the $n$-dependent version of Proposition \ref{prop-sec-mom-clt}(c)
in Proposition \ref{prop-sec-mom-clt}(d),
uniformly in $n$ and for each $i\in \BP^{\sss(j)}(s_n)$ and $k_n=o(\sqrt{\log{n}})$,
	\eqn{
	\e^{-\alphan t_n}
	\expec\Big[|\BP_{\sss \leq k_{t_n}(x)-k_n}^{\sss(i,j,\vec{m})}[t_n,t_n+s)|\Big]\ra A\Phi(x)\FR(s).	
	}
As a result, when $\tn+t-s_n-R_i\convp \infty$ and $G^{\sss(j)}_i=\op(\sqrt{\log{n}})$,
	\eqan{
	&\e^{-\alphan (\tn+t-s_n-R_i)}\expec\Big[|\BP^{\sss(i,j,\vec{m})}_{\sss \leq k_{t_n}(x)}[\tn+t-s_n-R_i, \tn+t+s-s_n-R_i)|\mid \FF_{s_n}\Big]\\
	&\qquad=\e^{-\alphan (\tn+t-s_n-R_i)}\expec\Big[|\BP^{\sss(i,j,\vec{m})}_{\sss \leq k_{t_n}(x)}[\tn+t-s_n-R_i, \tn+t+s-s_n-R_i)|\mid R_i,\tn\Big]\nn\\
	&\qquad=A\Phi(x)\FR(s)[1+\op(1)].\nn
	}
This yields that, when $G_i^{\sss(j)}=\op(\sqrt{\log{n}})$ 
(which happens \whp{} when $s_n$ is sufficiently small),
	\eqan{
	\label{SMM-mean}
	&\e^{-\alphan t_n}\sum_{i\in \BP^{\sss(j)}(s_n)}
	\expec\Big[\sum_{k=1}^{k_{t_n}(x)-G_i^{\sss(j)}} 
	|\BP_k^{\sss(i,j, \vec{m})}[\tn+t-s_n-R_i,\tn+t+s-s_n-R_i)|
	\mid \FF_{s_n}\Big]\\
	&\qquad =A\e^{\alpha t}\Phi(x)\FR(s) [1+\op(1)]\sum_{i\in \BP^{\sss(j)}(s_n)}\e^{\alphan(\tn-t_n-s_n-R_i)}\nn\\
	&\qquad \convp \e^{\alpha t} \Phi(x)\FR(s)\sqrt{\WW^{\sss(j)}/\WW^{\sss(3-j)}},\nn
	}
since $\e^{\alphan(t_n-\tn)}= \sqrt{\WW^{\sss(j)}_{s_n}\WW^{\sss(3-j)}_{s_n}}
\convp  \sqrt{\WW^{\sss(j)}\WW^{\sss(3-j)}}$, whereas by \eqref{Ms-convp}, $\sum_{i\in \BP^{\sss(j)}(s_n)}
	\e^{-\alpha_n(s_n+ R_i)}\convp \WW^{\sss(j)}/A$.

\paragraph{\ch{A conditional second moment method: second moment.}}
We next bound, conditionally on $\FF_{s_n}$, the variance of the sum on the right-hand side 
of \eqref{sum-trunc}. 
By conditional independence of $(\BP^{\sss(i,j)})_{i\geq 1}$,
	\eqan{
	&\e^{-2\alphan t_n}\Var\Big(\sum_{i\in \BP^{\sss(j)}(s_n)} \sum_{k=1}^{k_{t_n}(x)-G^{\sss(j)}_i}
	|\BP_k^{\sss(i,j, \vec{m})}[\tn+t-s_n-R_i,\tn+t+s-s_n-R_i)|
	\mid \FF_{s_n}\Big)\\
	&\qquad =\e^{-2\alphan t_n}
	\sum_{i\in \BP^{\sss(j)}(s_n)} \Var\Big(\sum_{k=1}^{k_{t_n}(x)-G^{\sss(j)}_i}|\BP_k^{\sss(i,j,\vec{m})}[\tn+t-s_n-R_i,\tn+t+s-s_n-R_i)|\mid \FF_{s_n}\Big).\nn
	}
Since adding individuals and enlarging the time-frame does not reduce the variance, 	
	\eqan{
	&\Var\Big(\sum_{k=1}^{k_{t_n}(x)-G^{\sss(j)}_i}|\BP_k^{\sss(i,j,\vec{m})}[\tn+t-s_n-R_i,\tn+t+s-s_n-R_i)|\mid \FF_{s_n}\Big)\\
	&\qquad \leq
	\expec\Big[|\BP^{\sss(i,j,\vec{m})}(\tn-s_n-R_i)|^2\mid \FF_{s_n}\Big]
	=\expec\Big[|\BP^{\sss(i,j,\vec{m})}(\tn-s_n-R_i)|^2\mid R_i,\tn\Big]\nn.
	}
By the $n$-dependent version of Proposition \ref{prop-sec-mom-clt}(b)
in Proposition \ref{prop-sec-mom-clt}(d), uniformly in $n$ and $t$,
	\eqn{
	\e^{-2\alphan t}\expec\big[|\BP^{\sss(i,j,\vec{m})}(t)|^2\big]\leq C\Kn.
	} 
As a result,
	\eqan{
	\label{SMM-var}
	&\e^{-2\alphan t_n}\Var\Big(\sum_{i\in \BP^{\sss(j)}(s_n)} \sum_{k=1}^{k_{t_n}(x)-G^{\sss(j)}_i}
	|\BP_k^{\sss(i,j,\vec{m})}[\tn+t-s_n-R_i,\tn+t+s-s_n-R_i)|
	\mid \FF_{s_n}\Big)\\
	&\qquad \leq \e^{-2\alphan t_n}
	\sum_{i\in \BP^{\sss(j)}(s_n)} \expec\Big[|\BP^{\sss(i,j,\vec{m})}(\tn-s_n-R_i)|^2\mid R_i,\tn\Big]\nn\\
	&\qquad  \leq C\Kn\e^{-2\alphan t_n}\sum_{i\in \BP^{\sss(j)}(s_n)} 	\e^{2\alphan (\tn-s_n-R_i)}\nn\\
	&\qquad =C\Kn\e^{2\alphan (\tn-t_n)} \e^{-2\alphan s_n}
	\sum_{i\in \BP^{\sss(j)}(s_n)} \e^{-2\alphan R_i}=\Op(1)\Kn\e^{-2\alphan s_n}
	\sum_{i\in \BP^{\sss(j)}(s_n)} \e^{-2\alphan R_i},\nn
	}
since $\e^{2\alphan (\tn-t_n)} =\Op(1)$. We can bound this further as
	\eqn{
	\Kn\e^{-2\alphan s_n}
	\sum_{i\in \BP^{\sss(j)}(s_n)} \e^{-2\alphan R_i}
	\leq \Kn\e^{-\alphan s_n}
	\Big(\e^{-\alphan s_n} \sum_{i\in \BP^{\sss(j)}(s_n)} 
	\e^{-\alphan R_i}\Big)=\Op(1) \Kn \e^{-\alphan s_n}
	=\op(1),
	}
precisely when $\Kn\e^{-\alphan s_n}=o(1)$.
By \eqref{SMM-mean} and \eqref{SMM-var},
the sum on the right-hand side of \eqref{sum-trunc} is, conditionally on $\FF_{s_n}$,
concentrated around its asymptotic conditional mean given in \eqref{SMM-mean}.
As a result, \eqref{aim-clt} follows.
This completes the proof of Proposition \ref{prop-CLT-stable-age}(a).
\hfill\qed
	
In order to prove Proposition \ref{prop-CLT-stable-age}(b), we need to 
investigate the asymptotics of the sum $\sum_{i=1}^{m} X^{\star}_i,$
where $m=|\SWG_{\sss\leq k_{\tn}(y)}^{\sss(j)}[\tn+t,\tn+t+s_2)|\convp \infty$ on the
event that $\WW_{s_n}^{\sss(1)}\WW_{s_n}^{\sss(2)}>0$, and 
$(X^{\star}_i)_{i\geq 1}$ are $X^{\star}_i=d_{V_i}-1$ with $(V_i)_{i\geq 1}$ the 
size-biased reordering of $(d_i)_{i\in [n]\setminus {\sf S}_m}$,
where ${\sf S}_m$ is the set of vertices found in 
$\SWG(\tn+t)$. We will prove that, conditionally on ${\mathcal F}_{\tn+t}$,
	\eqn{
	\label{conv-size-biased-sum}
	\frac{1}{m\nu_n}\sum_{i=1}^{m} X^{\star}_i\convp 1,
	}
and then the proof of Proposition \ref{prop-CLT-stable-age}(b) follows from 
the proof of Proposition \ref{prop-CLT-stable-age}(a) and the fact that $\nu_n\ra \nu$.

By the proof of Proposition \ref{prop-no-late-miscouplings-high-degree},
$|{\sf S}_m|\leq \mnup$, whp, where $\mnup=\sqrt{n}(\log{n})^{1/4}$ 
as in \eqref{mn-choices}. 
As a result the sequence $(d_i)_{i\in [n]\setminus {\sf S}_m}$ satisfies 
Condition \ref{cond-degrees-regcond}
when $(d_i)_{i\in[n]}$ does. Hence Lemma \ref{lem-mom-size-biased} holds with
$B_i$ replaced by $X_i^{\star}$, so that in  particular, from the Markov inequality, 
conditionally on ${\mathcal F}_{\tn+t}$,
	\eqn{
	\label{grotKnnul}
	\frac{1}{m}\sum_{i=1}^{m} X^{\star}_i\indic{X_i^{\star}>\Kn}\convp 0.
	}
We use a conditional second moment method on 
$\sum_{i=1}^{m} X^{\star}_i\indic{X_i^{\star}\leq \Kn}$,
conditionally on ${\mathcal F}_{\tn+t}$. 
By \eqref{size-biased-first} in 
Lemma \ref{lem-mom-size-biased},
	\eqn{
	\label{cond-sec-mom}
	\expec\Big[\sum_{i=1}^{m} X^{\star}_i\indic{X_i^{\star}\leq \Kn}\mid {\mathcal F}_{\tn+t}\Big]=m\nu_n(1+\op(1)).
	}
This gives the asymptotics of the first conditional moment of $\sum_{i=1}^{m} X^{\star}_i\indic{X_i^{\star}\leq \Kn}$. For the second moment,
we start by bounding the covariances. We note that, for $1\leq i<j\leq m$,
	\eqan{
	&{\rm Cov}\Big(X^{\star}_i\indic{X_i^{\star}\leq \Kn}, X^{\star}_j\indic{X_j^{\star}\leq \Kn}\mid {\mathcal F}_{\tn+t}\Big)\\
	&\quad =\expec\Big[X^{\star}_i\indic{X_i^{\star}\leq \Kn}
	\Big(\expec[X^{\star}_j\indic{X_j^{\star}\leq \Kn}\mid {\mathcal F}_{\tn+t}, X^{\star}_1,\ldots, X^{\star}_i]-\expec[X^{\star}_j\indic{X_j^{\star}\leq \Kn}\mid {\mathcal F}_{\tn+t}]\Big)\mid {\mathcal F}_{\tn+t}\Big].\nn
	}
By \eqref{size-biased-first} in Lemma \ref{lem-mom-size-biased}, as well as the fact that 
$i\leq \mnup=o(n)$, 
	\eqan{
	\expec[X^{\star}_j\indic{X_j^{\star}\leq \Kn}\mid {\mathcal F}_{\tn+t}, X^{\star}_1,\ldots, X^{\star}_i]-\expec[X^{\star}_j\indic{X_j^{\star}\leq \Kn}\mid {\mathcal F}_{\tn+t}]=\op(1),
	}
so that also
	\eqn{
	{\rm Cov}\Big(X^{\star}_i\indic{X_i^{\star}\leq \Kn}, 
	X^{\star}_j\indic{X_j^{\star}\leq \Kn}\mid {\mathcal F}_{\tn+t}\Big)=\op(1).
	}
Further, a trivial bound on the second moment together with
\eqref{size-biased-first} in Lemma \ref{lem-mom-size-biased} yields that	
	\eqn{
	\Var\Big(X^{\star}_i\indic{X_i^{\star}\leq \Kn}\mid {\mathcal F}_{\tn+t}\Big)
	\leq \Kn \expec[X^{\star}_i\mid {\mathcal F}_{\tn+t}]
	=\Kn\nu_n(1+\op(1)).
	}
As a result, whenever $\Kn m=o(m^2)$,
	\eqn{
	\Var\Big(\sum_{i=1}^{m} X^{\star}_i\indic{X_i^{\star}\leq\Kn}\mid {\mathcal F}_{\tn+t}\Big)=\op(m^2),
	}
which together with \eqref{cond-sec-mom} proves that, conditionally on ${\mathcal F}_{\tn+t}$, 
	\eqn{
	\frac{1}{m\nu_n}\sum_{i=1}^{m} X^{\star}_i\indic{X_i^{\star}\leq \Kn}\convp 1.
	}
Together with \eqref{grotKnnul}, this proves \eqref{conv-size-biased-sum},
as required.
\hfill\qed


\section{Extensions to other random graphs}
\label{sec-pfs-relat-RGs}

\paragraph{Proof of Theorem \ref{thm-unif-RGs}.} Let $\UGnd$ be a uniform
random graph with degree sequence $\bfd$. By \cite{Boll80b} (see also \cite{Boll01}),
we have that the law of $\UGnd$ is the same as that of $\CMnd$ conditioned
on  being simple, i.e., \ch{for every sequence of events $\HH_n$ defined on graphs with vertex set  $[n]$},
    \eqn{
    \prob(\UGnd\in \ch{\HH_n})=\prob(\CMnd\in \ch{\HH_n}\mid \CMnd \text{ simple})
    =\frac{\prob(\CMnd\in \ch{\HH_n},\CMnd \text{ simple})}{\prob(\CMnd \text{ simple})}.
    }
By \eqref{prob-simple}, it suffices to investigate $\prob(\CMnd\in \ch{\HH_n},\CMnd \text{ simple})$.
We take
    \eqn{
    \label{E-event}
    \ch{\HH_n}=\Big\{\frac{\Hn - \gamma_n\log{n}}{\sqrt{\beta\log{n}}}\leq x,
    \Wn - \frac{1}{\alphan} \log{n}\leq y\Big\},
    }
where $\Wn$ and $\Hn$ are the hopcount
and weight of the optimal path between two uniformly selected
vertices conditioned on being connected.

By the results in Section \ref{sec-pfs}, and with $\tilde{t}_n=t_n+\Bn$, \ch{where $\Bn=\log\log\log{n}$ is defined in 
\eqref{Bn-Cn-def}},
whp, we have found the minimal weight path before time
$\tilde{t}_n$. The probability that we have found
a self-loop or multiple edge at time $\tilde{t}_n$
is negligible, since, \ch{by that time
we have found of order $\mnup=\sqrt{n}(\log n)^{1/4}$ vertices and paired
of order $\mnup$ edges, see Lemma \ref{lem-particles-born}.}
Let $\tilde{d}_i(\tilde{t}_n)$ denote the number of unpaired
half-edges incident to vertex $i$ at time $\tilde{t}_n$.
Since $\CMnd$ is created by matching the half-edges
uniformly at random, in order the create $\CMnd$ after time $\tilde{t}_n$, we need to match the
half-edges corresponding to $(\tilde{d}_i(\tilde{t}_n))_{i\in[n]}$.
This corresponds to the configuration model on $[n]$ with degrees
$(\tilde{d}_i(\tilde{t}_n))_{i\in[n]}$. \ch{Since 
we have found of order $\mnup=\sqrt{n}(\log n)^{1/4}$ }vertices and paired
of order $\mnup$ edges at time $\tilde{t}_n$, when
$\bfd$ satisfies Condition \ref{cond-degrees-regcond},
then so does $(\tilde{d}_i(\tilde{t}_n))_{i\in[n]}$ with the \emph{same} limiting
degree distribution $D$. As result, the probability that
the configuration model on $[n]$ with degrees
$(\tilde{d}_i(\tilde{t}_n))_{i\in[n]}$ is simple is asymptotically equal to
$\e^{-\nu/2-\nu^2/4}(1+o(1))$, and we obtain that the event that
$\CMnd$ is simple is asymptotically independent of the event
$\ch{\HH_n}$ in \eqref{E-event}. Therefore, Theorem \ref{thm-unif-RGs}
follows from Theorems \ref{thm-main-first}-\ref{thm-main}.
\hfill \qed

\paragraph{Proof of Theorem \ref{thm-IRGs}.}
By Janson \cite{Jans08a}, when $W_n\convd W$ and
$\expec[W_n^2]\rightarrow \expec[W^2]$, the inhomogeneous random graphs
with edge probabilities in \eqref{pij-def},
\eqref{pij-GRG} or \eqref{pij-CL}
are asymptotically equivalent, so it suffices to prove the
claim for the \emph{generalized random graph} for which
$p_{ij}=w_iw_j/(\ell_n+w_iw_j)$.
As explained in Section \ref{sec-rel-models}, conditionally on the degrees in
the generalized random graph being equal to $\bfd$, the distribution
of the resulting random graph is uniform over all random graphs with
these degrees. Therefore,  Theorem \ref{thm-IRGs} follows from
Theorem \ref{thm-unif-RGs} if we prove that Condition \ref{cond-degrees-regcond}
follows from the statements that
$W_n\convd W$, $\expec[W_n]\rightarrow \expec[W],
\expec[W_n^2]\rightarrow \expec[W^2]$ and \ch{$\lim_n \expec[W_n^2\log{(W_n/\Kn)_+}]=0$}.
We denote by $\bfd=(d_i)_{i\in[n]}$ the degree sequence in the generalized random graph,
and note that $\bfd$ now is a \emph{random} sequence. We work conditionally on
$\bfd$, and let $\prob_n,\expec_n$ denote the conditional probability and
expectation given $\bfd$. Then, we prove that $\prob_n(D_n=k)\convp \prob(D=k),
\expec_n[D_n]\convp \expec[D],
\expec_n[D_n^2]\convp \expec[D^2]$ and $\expec_n[D_n^2 \log(D_n/\Kn)_+]\convp 0$
for every $\Kn\ra \infty$.

We let $D_n=d_{V}$, where $V\in[n]$ is a uniformly chosen vertex.
First, by \eqref{degree-conv-IRG}, $D_n\convd D$, \ch{where $D$ is a Poisson random variable with random intensity $W$}. 
Further,
    \eqn{
    \expec_n[D_n]=\frac{1}{n} \sum_{i\in [n]} d_i,
    \qquad
    \expec_n[D_n^2]=\frac{1}{n} \sum_{i\in [n]} d_i^2,
    }
where $d_i=\sum_{j\in [n], j\neq i} I_{ij}$ and $I_{ij}$ are independent Bernoulli variables
with parameter $p_{ij}=w_iw_j/(\ell_n+w_iw_j)$. It tedious, but not difficult, to show that
the above sums are concentrated around their means, for example by computing their means
and variances and showing that the variances are of smaller order than their means squared.
We omit the details.

\ch{In order to show that $\expec_n[D^2_n \log(D_n/\Kn)_+]=\op(1)$}, we note that
    \eqn{
    \expec_n[D_n^2 \log(D_n/\Kn)_+]
    =\frac{1}{n} \sum_{i\in [n]} d_i^2\log{(d_i/\Kn)_+}.
    }
As before, $d_i=\sum_{j\in [n], j\neq i} I_{ij}$ and $I_{ij}$ are independent Bernoulli variables
with parameter $p_{ij}=w_iw_j/(\ell_n+w_iw_j)$. By standard Chernoff bounds, 
there exists a constant $a>0$ such that, for every $\lambda>2$,
    \eqn{
    \label{chernovforbin}
    \prob(d_i\geq \lambda\expec[d_i])\leq \e^{-a \lambda \expec[d_i]}.
    }
Here,
    \eqn{
    \label{bounddbyw}
    \expec[d_i]=\sum_{j\neq i} w_iw_j/(\ell_n+w_iw_j)\in w_i(\frac{\ell_n-\sqrt{n}}{\ell_n(1+o(1))}, 1),
    }
since $\max_{i}w_i=o(\sqrt{n})$.
\ch{As a result, 
    \eqn{
	\label{XlogX-GRGa}
    \expec_n[D_n^2 \log(D_n/\Kn)_+]
    \leq \frac{4}{n} \sum_{i\in [n]} w_i^2\log{(2w_i/\Kn)_+}
    +\frac1n\sum_{i\in [n]} \indic{d_i\geq 2w_i} d_i^2\log{(d_i/\Kn)_+}.
    }
    }
The first term vanishes by the fact that \ch{$\lim_{n\ra \infty}\expec[W^2_n\log{(W_n/\Kn)_+}]=0$}.
The second term can be split as \ch{
    \eqan{
	\label{XlogX-GRGb}
    \frac1n \sum_{i\in [n]} \indic{d_i\geq 2w_i} d_i^2\log{(d_i/\Kn)_+}
    & \leq \frac1n \sum_{k=1}^{\infty} \sum_{i\in [n]} \indic{d_i\in [2^kw_i, 2^{k+1}w_i)} d_i^2\log{(d_i/\Kn)_+}\\
    &\leq \frac1n\sum_{i\in [n]} \sum_{k=1}^{\infty} 4^{k+1} \indic{d_i\geq 2^{k}w_i} w_i^2\log{(2^{k+1}w_i/\Kn)_+}.
    \nn
    }
    }
\ch{By \eqref{chernovforbin}-\eqref{bounddbyw}} with $\lambda=2^k$, 
the mean of the above random variable vanishes, which, by Markov's inequality,
implies that it converges to zero in probability.
\hfill\qed

\paragraph{Acknowledgements.}
The research of SB is supported by NSF-DMS grant 1105581
and would like to thank Eurandom and TU/e for travel support where part of this work was done.
The work of RvdH was supported in part by Netherlands
Organisation for Scientific Research (NWO). SB and RvdH thank Cornell Unversity for 
its hospitality during the Cornell Probability Summer School in July 2012, where
this work was completed.

\bibliographystyle{plain}

\end{document}